\RequirePackage{etex}

\documentclass[12pt,psamsfonts]{amsart}
\usepackage{amsfonts,pstricks}
\usepackage{amsmath, amsthm, amssymb,  mathrsfs, epsfig}
\usepackage{dcpic, pictex, pinlabel}
\usepackage{hyperref,url}
\usepackage{graphicx}
\usepackage{color}
\usepackage{ifpdf}
\usepackage{tikz}
\usepackage{xypic}
\usepackage{extarrows}
\usepackage{multirow}
\usepackage{bm}

\begin{document}
\newcommand{\Z}{\mathbb{Z}}
\newcommand{\M}{\mathcal{M}}
\newcommand{\B}{\mathcal{B}}
\newcommand{\C}{\mathcal{C}}
\newcommand{\Sym}{\mathcal{S}}
\newcommand{\A}{\mathcal{A}}
\newcommand{\Q}{\mathcal{Q}}
\def\hom#1{\phi_{#1}}

\def\Ap{\mathcal{A}_{n}^{+}}
\def\Am{\mathcal{A}_{n}^{-}}
\def\An{\mathcal{A}_{n}}
\def\Apm{\mathcal{A}_{n}^{\pm}}

\def\Phim#1{\Phi_{#1}^{-}}
\def\Phip#1{\Phi_{#1}^{+}}
\def\Phipm#1{\Phi_{#1}^{\pm}}

\def\Phiml#1{\Phi_{#1}^{-L}}
\def\Phimr#1{\Phi_{#1}^{-R}}
\def\Phipl#1{\Phi_{#1}^{+L}}
\def\Phipr#1{\Phi_{#1}^{+R}}

\def \Ep{\epsilon^{+}}
\def \Em{\epsilon^{-}}
\def \Epm{\epsilon^{\pm}}

\def \Var{\Gamma}

\newcommand{\Lens}{S^1 \times S^2}
\theoremstyle{definition}
\newtheorem{axiom}{Axiom}
\newtheorem{thm}{Theorem}
\newtheorem{conjecture}{Conjecture}
\newtheorem{lem}{Lemma}
\newtheorem{example}{Example}
\newtheorem{cor}{Corollary}
\newtheorem{prop}{Proposition}
\newtheorem{rem}{Remark}
\newtheorem{definition}{Definition}
\newtheorem{measurement}{Measurement}
\newtheorem{ancilla}{Ancilla}

\numberwithin{equation}{section} \makeatletter
\renewenvironment{proof}[1][\proofname]{\par
    \pushQED{\qed}%
    \normalfont \topsep6\p@\@plus6\p@ \labelsep1em\relax
    \trivlist
    \item[\hskip\labelsep\indent
        \bfseries #1]\ignorespaces
}{%
    \popQED\endtrivlist\@endpefalse
} \makeatother
\renewcommand{\proofname}{Proof}

\title{Framed Cord Algebra Invariant of Knots in $\Lens$}
\author{Shawn X. Cui$^{1}$ and Zhenghan Wang$^{1,2}$}
\address{$^1$Department of Mathematics\\University of California\\Santa Barbara, CA 93106}
\email{xingshan@math.ucsb.edu, zhenghwa@math.ucsb.edu}
\address{$^2$Microsoft Research, Station Q\\ University of California\\ Santa Barbara, CA 93106}
\email{zhenghwa@microsoft.com}

\thanks{The authors are partially supported by NSF DMS 1108736. Both authors thank Lenny Ng for valuable comments on the paper.  He pointed out that the cord algebras of local knots are infinitely generated, and that the ends of the cords should stay on the framing.}


\maketitle

\begin{abstract}

We generalize Ng's two-variable algebraic/combinatorial $0$-th framed knot contact homology for framed oriented knots in $S^3$ to knots in $\Lens$, and prove that the resulting knot invariant is the same as the framed cord algebra of knots.  Actually, our cord algebra has an extra variable, which potentially corresponds to the third variable in Ng's three-variable knot contact homology.  Our main tool is Lin's generalization of the Markov theorem for braids in $S^3$ to braids in $\Lens$.  We conjecture that our framed cord algebras are always finitely generated for non-local knots.

\end{abstract}

\section{Introduction}
\label{sec:introduction}

The dream of finding new higher categorical quantum invariants of smooth $4$-manifolds that can distinguish smooth structures beyond Donaldson/Seiberg-Witten/Heegaard-Floer theory is largely unrealized, despite the spectacular success for new invariants in $3$-dimensions and recent progress in higher category theory.  A potentially new quantum invariant would be to promote the relative knot contact homology of knots in $S^3$ in \cite{ng2008framed} to a $(3+1)$-TQFT-type theory (presumably the $0$-th part of the BRST cohomology of a topological string theory).  One lesson from $(2+1)$-dimensions is the emergence of powerful diagrammatical techniques as exemplified by the Kauffman bracket definition of the Jones polynomial, and the subsequently elementary formulation of Turaev-Viro and Reshetikhin-Turaev $(2+1)$-TQFTs.  We see a striking parallel between the cord algebra invariant and the Jones polynomial.

In \cite{ng2008framed}, the $0$-th part of the relative knot contact homology in $S^3$ is interpreted using cords and skein relations---the main ingredients of diagrammatical techniques in $(2+1)$-dimensions, analogous to the reformulation of the Jones polynomial of knots from von Neumann algebra using knot diagrams and the Kauffman bracket.  Taking the elementary framed cord algebra invariant of knots in general $3$-manifolds $M$ as the main object of interest, we will follow the diagrammatical approach to constructing $(2+1)$-TQFTs such as the Turaev-Viro and Reshetikhin-Turaev TQFTs.  As a first step, we generalize Ng's two-variable combinatorial/algebraic $0$-th framed knot contact homology for framed oriented knots in $S^3$ to knots in $\Lens$, and prove that the resulting knot invariant is the same as the framed cord algebra of knots.  Actually, our cord algebra has an extra variable, which potentially corresponds to the third variable in Ng's three-variable knot contact homology \cite{ng2011combinatorial}.

It is conjectured in \cite{ng2005knot} that the cord algebra invariant of knots in a general $3$-manifold $M$ is the $0$-th relative knot contact homology.  We do not prove this conjecture and will not use any knot contact homology theory.  Instead we provide an algebraic version of this conjectured $0$-th knot contact homology for knots in $\Lens$ following \cite{ng2008framed} and regard our algebraic definition of the cord algebra as an effective method to calculate the topologically defined cord algebra invariant of knots.  Our long term goal is to understand the higher categories underlying this algebraic formulation with an eye towards to a diagram construction of a $(3+1)$-TQFT-type theory.

A second reason for our interest in the framed cord algebra invariant of knots is the conjectured relation between the augmentation polynomial and the Homfly polynomial of knots.  A well-known question since the discovery of the Jones polynomial is how to place the Jones polynomial within classical topology (since knots are determined by their complements, so any knot invariant is determined by the homeomorphism type of the knot complement).  The cord algebra of a knot is basically within classical topology, so the establishment of the conjectured relation between the augmentation polynomial and the Homfly polynomial is one answer to an old question.

To generalize the algebraic $0$-th knot contact homology in \cite{ng2008framed} from $S^3$ to $\Lens$, we use Lin's generalization of the Markov theorem for braids in $S^3$ to braids in $\Lens$ \cite{Lin} developed for defining a Jones polynomial of knots in $\Lens$.\footnote{This generalization, eventually rendered unnecessary for the intended application by Witten's work, finds a similar application in our work.  We dedicate our work to X.-S. Lin---an important vanguard in quantum knot theory.}

The rest of the paper is organized as follows. In Section \ref{subsec:markov}, we introduce the Markov theorem for knots in $\Lens$, which are represented by the closure of elements in $\C_n$, the Artin group with Dynkin diagram $B_n$. In Section \ref{subsec:braid actions}, we give several actions of $\C_n$ on free algebras. We interpret these actions both algebraically and topologically. These actions will be the key ingredients to define the invariant $HC_0$ in Section \ref{subsec:def invariant}. In Section \ref{subsec:example} - Section \ref{subsec:invariace proof}, we compute some specific examples, demonstrate some useful propositions, and prove the invariance of $HC_0$ under Markov moves, respectively. Section \ref{subsec:symmetry} - Section \ref{subsec:augmentation} are devoted to prove several properties of the $HC_0$ invariant. We study two special classes of knots in $\Lens$, torus knots and local knots. Moreover, we derive a family of invariants, called augmentations, from $HC_0$. Finally, in Section \ref{sec:topology} we prove that the $HC_0$ invariant has a nice topological interpretation as the framed cord algebra defined in \cite{ng2008framed}.

The first author also created a Mathematica package for computer calculations of the $HC_0$ invariant and augmentation numbers. The program can be found at \cite{CuiPackage} and is partly motivated by Ng's computer package, which was used to compute various invariants derived from knot contact homology for knots in $S^3$. To run the program, one needs to install the non-commutative algebra package NCAlgebra$/$NCGB \cite{HeltonPackage}.

\section{Markov moves and actions of $\C_n$ on free algebras}
\label{sec:preliminary}

First we provide some background materials.  Links and knots in this paper are always framed and oriented.

\subsection{Markov moves in $\Lens$}
\label{subsec:markov}

In this subsection, we describe a theorem on Markov moves for links in $\Lens$. See \cite{Lin} for a more detailed discussion.

The classical braid group with $n$ strands, $\B_n$, is defined by the presentation $\langle \sigma_1, \cdots, \sigma_{n-1} | \sigma_{i}\sigma_{i+1}\sigma_{i} = \sigma_{i+1}\sigma_{i}\sigma_{i+1}, \sigma_{i}\sigma_{j} = \sigma_{j}\sigma_{i}, |i-j| \geq 2 \rangle.$ It is the Artin group with Dynkin diagram type $A_{n-1},$ and can also be viewed as the braid group on the $2$-disk $D^2 \subset \mathbb{R}^2.$

Any link in $S^3$ can be represented as the closure of some braid in the classical braid group. The Markov theorem states that two braids $B, B'$ give rise to the same link if and only if $B'$ can be obtained from $B$ by a finite sequence of the following operations or their inverses:

$1).$ change $B \in \B_n$ to one of its conjugates in $\B_n;$

$2).$ change $B \in \B_n$ to $B\sigma_n^{\pm 1} \in \B_{n+1}.$

The Markov theorem for links in $S^3$ is generalized to links in $\Lens$ in \cite{Lin} as follows.

Let $\C_{n}$ be the Artin group corresponding to the Dynkin diagram $B_n$ generated by $\alpha_0, \cdots, \alpha_{n-1},$ with the following generating relations:

$1). \;  \alpha_{i}\alpha_{j} = \alpha_{j}\alpha_{i}, \; |i-j| \geq 2$

$2). \; \alpha_{i}\alpha_{i+1}\alpha_{i} = \alpha_{i+1}\alpha_{i}\alpha_{i+1}, \; i \geq 1$

$3). \; \alpha_{0}\alpha_{1}\alpha_{0}\alpha_{1}=\alpha_{1}\alpha_{0}\alpha_{1}\alpha_{0}.$

Clearly, we have natural inclusions $\C_1 \subset \C_2 \subset \cdots \subset \C_n \subset \cdots $. We denote by $\Em$ these natural embeddings.

\begin{figure}[h!]
\centering
\def\svgwidth{10cm}
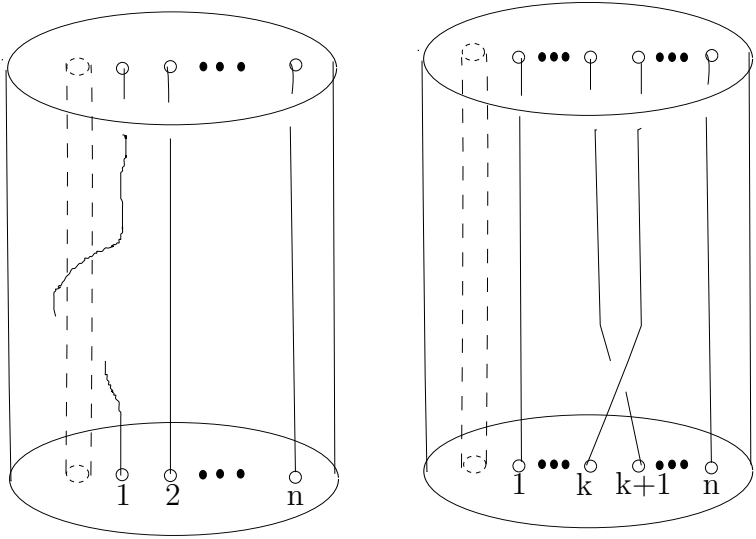
\caption{$\alpha_0$ \textrm{and} $\alpha_k, k \geq 1$}
 \label{alpha_0k}
\end{figure}

It is shown in \cite{crisp1999injective} that $\C_n$ is isomorphic to the braid group on the annulus $[0,1]\times S^1,$ or the $1$-punctured disk. Specifically, the isomorphism is illustrated in Figure \ref{alpha_0k}.

Simply treating $\{\textrm{puncture}\} \times [0,1]$ as the first strand of the new braid, we can regard a braid on the $1$-punctured disk as a braid on the disk. Thus we have an embedding of $\C_n$ into $\B_{n+1}.$ Denote the generators of $\B_{n+1}$ by $\sigma_0,\sigma_1, \cdots, \sigma_{n-1}.$ Then the embedding from $\C_n$ to $\B_{n+1}$ is given by by the following map:
$$
\C_n \longrightarrow \B_{n+1}, \quad \alpha_0 \longmapsto \sigma_0^2, \quad \alpha_i \longmapsto \sigma_i, i\geq 1.
$$

From now on, we will identify $\C_n$ with its image in $B_{n+1}$, which is the subgroup consisting of the braids that fix the first puncture.

The correspondence between braids on the annulus and links in $\Lens$ is obtained via open book decompositions.

Consider the standard open book decomposition of $S^3$ with an unknot $J$ as the binding. Let $K$ be another unknot which is a closed braid with respect to the braid axis $J$. Then

$$
M = \overline{S^3 \setminus (J \times D^2  \cup K \times D^2)}
$$
is a fibration over $S^1$ whose fibre is an annulus $[0,1] \times S^1$. $\Lens$ is obtained by a $0$-Dehn surgery along $K$. Thus $\Lens = M \sqcup_{f} D^2 \times S^1$, where $f$ is the gluing homeomorphism which maps the meridian of the solid torus to $K \times z_0, z_0 \in \partial D^2.$   Let $K^{*}$ be the image of ${0} \times S^1$ under $f$ in $S^1 \times S^2$, where ${0} \times S^1$ is the core of the solid torus.  We call $K^{*}$ the dual knot of $K$.  Then the fibration on M extends to an open book decomposition on $\Lens$ with the binding $J \cup K^{*}$. Note that $\Lens \setminus (J \cup K^{*})$ is homeomorphic to the product of the annulus with $S^1$. It's not hard to see that any link in $\Lens$ can be isotoped into $\Lens \setminus (J \cup K^{*})$ transversal to each page, and thus becomes a braid on the annulus.

To state the Markov theorem, we need one more lemma.

Define a map $\Ep: \C_n \longrightarrow \C_{n+1},$

\begin{equation}
\Ep(\alpha_i) =
\begin{cases}
\alpha_1\alpha_0\alpha_1  & i = 0 \\
\alpha_{i+1}              & i \geq 1 \\
\end{cases}
\end{equation}

\begin{figure}[h!]
\centering
\def\svgwidth{4cm}
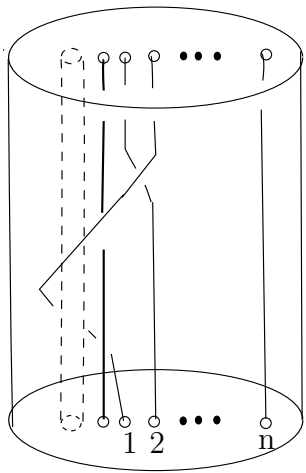
\caption{$\epsilon^{+}(\alpha_1\alpha_0)$}
 \label{epsilonplus}
\end{figure}

The map $\Ep$ has a nice geometrical interpretation if we view $\C_n$ as the braid group on the annulus. The map simply inserts a straight strand right next to the line $\{\textrm{puncture}\} \times [0,1]$. See Figure \ref{epsilonplus}.

Note that the newly inserted line will be labeled by $1$, and the other strands' labels will be shifted up by $1$.

\begin{lem}\cite{Lin}
The map $\Ep$ is an injective group homomorphism.
\begin{proof}
From the geometrical interpretation of the map, it should be clear that it is an injective group homomorphism. For a rigorous algebraic proof, see \cite{Lin}.
\end{proof}
\end{lem}

\begin{rem}
Now there are two embeddings of $\C_n$ into $\C_{n+1},$ namely the natural inclusion $\Em$ and the map $\Ep$. From the geometric point of view, $\Em$ is to place a strand on the far right of the braid, while $\Ep$ is to insert a strand right next to the line $\{\textrm{puncture}\} \times [0,1]$.
\end{rem}

Here is the statement of the Markov Theorem for links in $\Lens$.

\begin{thm} \label{markov move} \cite{Lin}
The closures of two braids $\beta, \beta' \in \cup_{n=1}^{\infty} \C_{n}$ give the same link in $\Lens$ if and only if  there is a finite sequence of braids, $\beta = \beta_0, \beta_1, \cdots, \beta_k = \beta'$, such that $\beta_{i+1}$ can be obtained from $\beta_i$ by one of the following operations or their inverses:

$1).$ change $\beta_i \in \C_n$ to one of its conjugates in $\C_n$;

$2).$ change $\beta_i \in \C_n$ to $\Em(\beta_i)\alpha_{n}^{\pm} \in \C_{n+1}$;

$3).$ change $\beta_i \in \C_n$ to $\Ep(\beta_i)\alpha_{1}^{\pm} \in \C_{n+1}$.
\end{thm}

\begin{rem} \label{S1S2 model}
Given a braid $\beta \in \C_n$, we can obtain the knot in $\Lens$ represented by $\beta$ as follows. Take a punctured disk $D' = D \setminus B_{\epsilon}(0)$, and let $X = D' \times [0,1]$. Draw the diagram of $\beta$ inside $X$. Then $\Lens$ is obtained by identifying the top and the bottom punctured disk and then gluing a solid torus to each torus boundary component. The gluing maps are given by sending the meridian of each solid torus to $z_0 \times S^1$ and $z_1 \times S^1$, respectively for some $z_0$ on the boundary of the puncture and $z_1$ on the outer boundary of $D'$. And the knot represented by $\beta$ is the image of the braid diagram in $\Lens$. See Figure \ref{ModelS1S2} for $\beta = \alpha_0\alpha_1$.

\begin{figure}[h!]
\centering
\def\svgwidth{5cm}
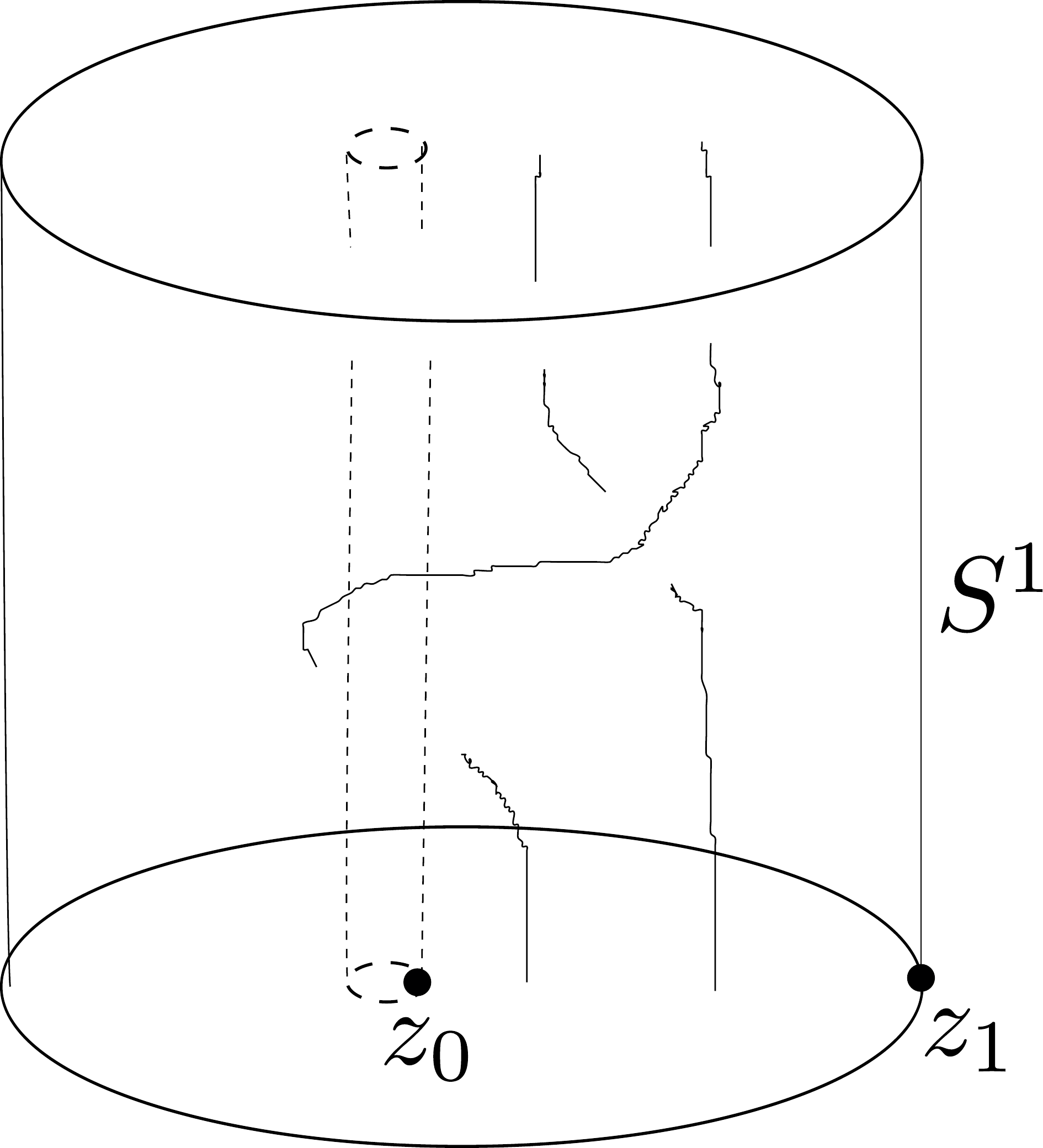
\caption{The closure of $\beta$ in $\Lens$}
 \label{ModelS1S2}
\end{figure}

\end{rem}

\subsection{Actions of $\C_n$ on free algebras}
\label{subsec:braid actions}

Throughout the paper, $R$ denotes the commutative ring $\Z[\lambda^{\pm},\mu^{\pm},\Var^{\pm}].$  We define several free non-commutative algebras over the ring $R$ as follows.

$\Ap := R \langle a_{ij}^{x}, 0 \leq i,j \leq n, x \in \Z\rangle / \langle a_{ii}^0 - (1+\mu)\Var, 0 \leq i \leq n \rangle$,

$\Am := R \langle a_{ij}^{x}, 1 \leq i,j \leq n+1, x \in \Z\rangle / \langle a_{ii}^0 - (1+\mu)\Var, 1 \leq i \leq n+1 \rangle$,

$\An := R \langle a_{ij}^{x}, 1 \leq i,j \leq n, x \in \Z\rangle / \langle a_{ii}^0 - (1+\mu)\Var, 1 \leq i \leq n \rangle.$

The algebra $\An$ can be embedded into $\Ap$ and $\Am$ in the most natural way. We will always identity $\An$ with its images in $\Ap$ and $\Am$.

Now we introduce an action of $\C_n$ on $\A_n$, and extend the action to the larger algebras $\Ap, \Am$. The action is first presented algebraically and then will be given a topological interpretation.

Recall that the generators $\C_n$ are denoted by $\alpha_0, \cdots, \alpha_{n-1}$, which satisfy the relation given in Section \ref{subsec:markov}. We define a group morphism $\Phi :  \C_n \longrightarrow \textrm{Aut}(\A_n)$ as follows.

For $1 \leq k \leq n-1,$

\begin{equation} \label{equ:Phi_k}
\Phi(\alpha_k)(a_{ij}^x) =
\begin{cases}
-a_{k+1,j}^x + \frac{1}{\Var\mu}a_{k+1,k}^0 a_{k,j}^x              &  i = k, j  \neq k,k+1 \\
-a_{k+1,k}^x + \frac{1}{\Var\mu}a_{k+1,k}^0 a_{k,k}^x                  &  i = k, j = k+1 \\
a_{k+1,k+1}^x - \frac{1}{\Var}a_{k+1,k}^x a_{k,k+1}^0 -             &           \\
\quad  \frac{1}{\Var\mu}a_{k+1,k}^0 a_{k,k+1}^x + \frac{1}{\Var^2\mu}a_{k+1,k}^0 a_{k,k}^x a_{k,k+1}^0               &  i = k, j =k \\
a_{k,j}^x                                                      &  i=k+1, j \neq k,k+1 \\
-a_{k,k+1}^x + \frac{1}{\Var}a_{k,k}^x a_{k,k+1}^0                  &  i = k+1, j = k \\
a_{k,k}^x                                                      &  i=k+1, j =k+1 \\
-a_{i,k+1}^x + \frac{1}{\Var}a_{i,k}^x a_{k,k+1}^0                  &  i \neq k, k+1, j = k \\
a_{i,k}^x                                                      &  i \neq k,k+1, j =k+1 \\
a_{i,j}^x                                                      &  i \neq k,k+1, j \neq k,k+1 \\
\end{cases}
\end{equation}

\begin{equation} \label{equ:Phi_0}
\Phi(\alpha_0)(a_{ij}^x) =
\begin{cases}
a_{1,1}^{x}                                                   & i = 1, j = 1 \\
-\mu a_{1,j}^{x-1} + \frac{1}{\Var}a_{1,1}^{x}a_{1,j}^{-1}           & i = 1, j \geq 2 \\
\frac{1}{\mu}(-a_{i,1}^{x+1} + \frac{1}{\Var}a_{i,1}^{1}a_{1,1}^{x})            & i \geq 2, j =1 \\
a_{i,j}^x - \frac{1}{\Var\mu}a_{i,1}^{x+1} a_{1,j}^{-1} -          &           \\
\quad  \frac{1}{\Var}a_{i,1}^1 a_{1,j}^{x-1} + \frac{1}{\Var^2\mu}a_{i,1}^1 a_{1,1}^x a_{1,j}^{-1}               &  i \geq 2, j \geq 2 \\
\end{cases}
\end{equation}

It is not hard, though tedious, to check that $\Phi$ is well defined, i.e. $\Phi(\alpha_i)$ satisfies the braid relations that define $\C_n.$

We extend the action of $\C_n$ to the algebra $\Ap$ by furthermore defining the action on $a_{0j}^x, a_{i0}^x, 0 \leq i,j \leq n.$  This extended action will be denoted by $\Phi^{+}.$

\begin{equation}
\Phi^{+}(\alpha_0)(a_{ij}^x) =
\begin{cases}
a_{0,0}^{x}                                                   & i = 0, j = 0 \\
\frac{1}{\mu} a_{0,1}^{x+1}                                                 & i = 0, j = 1 \\
-a_{0,j}^x + \frac{1}{\Var\mu}a_{0,1}^{x+1}a_{1,j}^{-1}             & i = 0, j \geq 2 \\
\mu a_{1,0}^{x-1}                                                 & i = 1, j = 0 \\
-a_{i,0}^x + \frac{1}{\Var}a_{i,1}^{1}a_{1,0}^{x-1}              & i \geq 2, j =0 \\
\end{cases}
\end{equation}

For $1 \leq k \leq n-1$, $\Phi^{+}(\alpha_k)(a_{ij}^x)$ are given by the same equation as \ref{equ:Phi_k}, except that now $i,j$ are allowed to be zero when they are not $k$ or $k+1$.

Similarly, the extended action of $\C_n$ on $\Am$ is defined by Equations \ref{equ:Phi_k}, \ref{equ:Phi_0} except that the range of $i,j$ now is from $1$ to $n+1$. We denote this action by $\Phi^{-}.$

Again, it can be checked $\Phi^{+}, \Phi^{-}$ are both well defined.

A few remarks are in order.

\begin{rem} \label{rem beta action}
1). From now on, for a braid $\beta \in C_n$, we will write $\Phi_{\beta}, \Phip{\beta}, \Phim{\beta}$ for $\Phi(\beta), \Phi^{+}(\beta), \Phi^{-}(\beta)$, respectively.

2). It's direct from the very definitions that $\Phi_{\beta} = (\Phip{\beta})_{|\An} =  (\Phim{\beta})_{|\An}.$ It's also clear that $\Phim{\beta} = \Phi_{\Em(\beta)}$ if we identify $\Am$ with $\A_{n+1}$ in the obvious way.

3). From the definition of $\C_n$ in Section \ref{subsec:markov}, it's easy to see that the subgroup generated by $\{\alpha_1, \cdots, \alpha_{n-1}\}$ is isomorphic to the classical braid group on $n$ strands. We denote this subgroup by $\B_n$. In Equation \ref{equ:Phi_k}, if we set $\Var=-1, \mu = 1,$ and $x = 0$, then $\Phi_{| \B_n}$ acting on $\Z\langle a_{ij}^0 \rangle$ is exactly the braid group action given in \cite{ng2005knotI}. So our braid group action is a generalization of Ng's in \cite{ng2005knotI}.

\end{rem}

The above actions will be less mysterious after we give a topological interpretation.

\begin{figure}[h!]
\centering
\def\svgwidth{5cm}
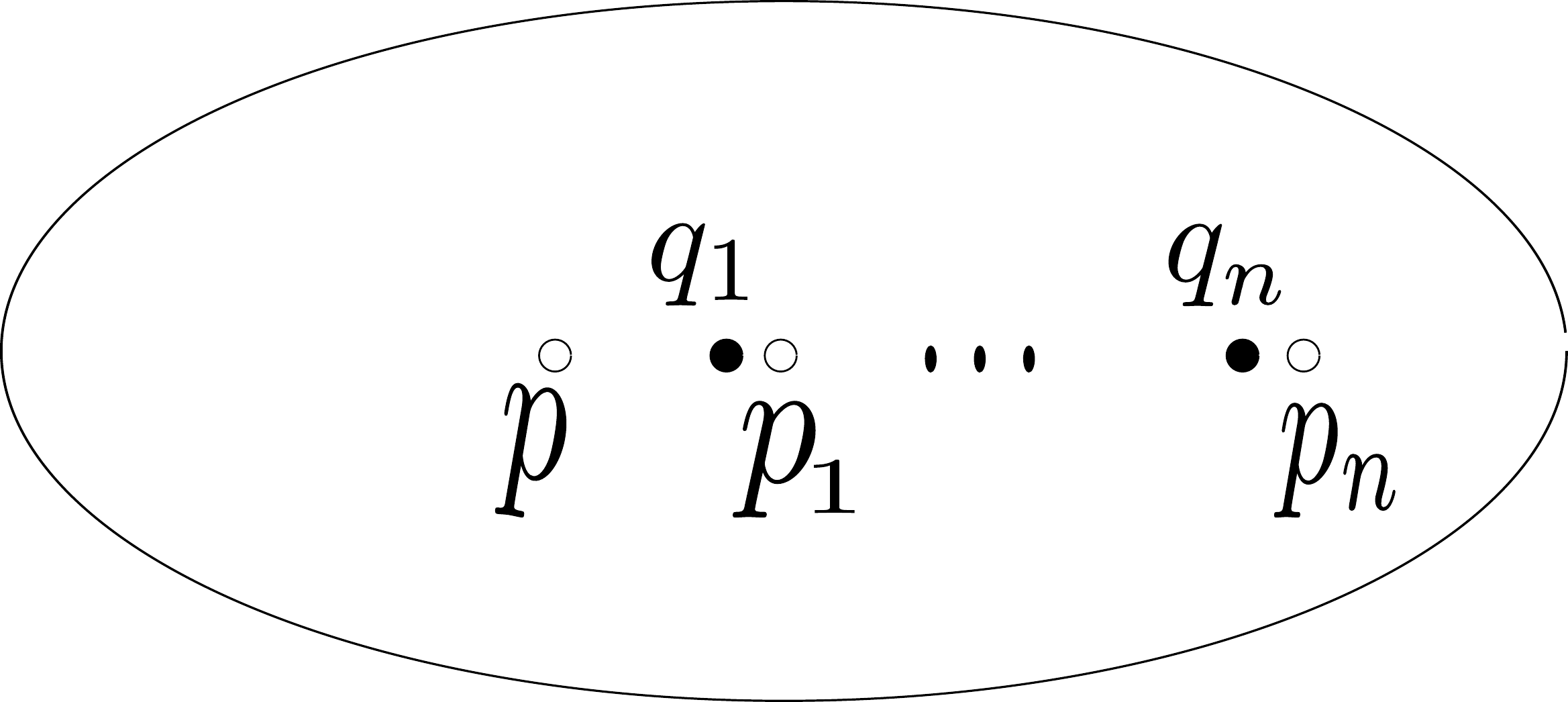
\caption{$D_n$}
\label{D_n}
\end{figure}

Let $D$ be the unit disk in the complex plane centered at the origin, $D_n$ be the punctured disk with $n+1$ punctures labeled, from left to right, by $p, p_1, \cdots, p_{n}$ and let $q_i = p_i - \epsilon, 1 \leq i \leq n$, $\epsilon > 0$ be $n$ points in $D_n$ close to the punctures. See Figure \ref{D_n}.

Let $Q_n = \{q_i, 1 \leq i \leq n\}$ and let $\Q_n = \{\gamma:[0,1] \longrightarrow D_n | \, \gamma  \\  \textrm{is continuous}, \gamma(0),\gamma(1) \in Q_n\} / \thicksim$. Here $\thicksim$ is the equivalence relation which means two curves $\gamma_1 \thicksim \gamma_2$ if and only if $\gamma_1$ and $\gamma_2$ are homotopic inside $D_n$ relative to their end points. So the curves are not allowed to pass through any of the punctures and their end points are fixed during the homotopy. Then $\Q_n$ is the  set of equivalence classes of such curves.

\begin{figure}[h!]
\def\svgwidth{11cm}
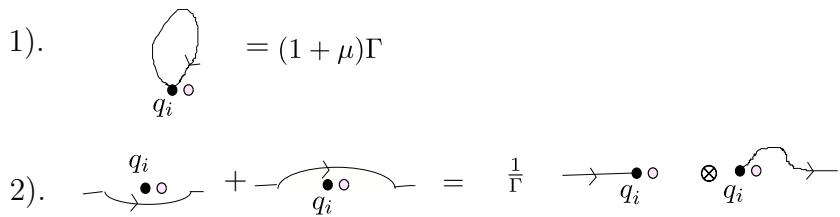
\caption{\textrm{skein relation}}
\label{skein12}
\end{figure}

Let $\tilde{\A}_{n}$ be the free non-commutative algebra over $R$ generated by elements of $\Q_n$ modulo the \lq\lq skein" relations shown in Figure \ref{skein12}. Note that $\otimes$ in Figure \ref{skein12}, and all other places of the paper,  means the multiplication. And the second relation, as well as other similar relations in the context, depicts some local neighborhood of the diagrams outside of which they all agree.

\begin{figure}[h!]
\centering
\def\svgwidth{11cm}
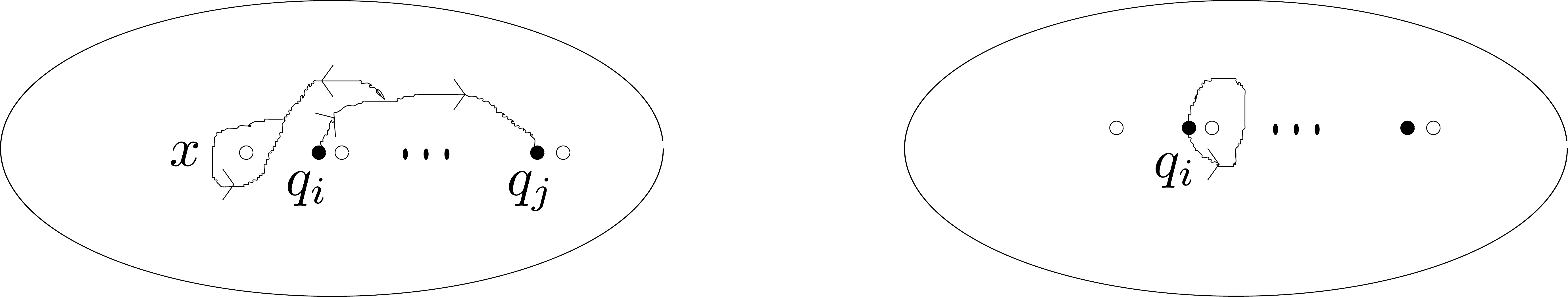
\caption{$\gamma_{ij}$ and $\gamma_i$}
 \label{gammaij}
\end{figure}

For $1 \leq i,j \leq n, x \in \Z$, let $\gamma_{ij}^x$ and $\gamma_i$ be the curves shown in Figure \ref{gammaij}, namely $\gamma_{ij}^x$ starts from $q_i$, winds around $p$ counter clock-wise $x$ times if $x \geq 0$, or clock-wise $-x$ times if $x < 0,$ and finally goes through the upper half disk to reach $q_j$. $\gamma_i$ is the curve that starts and ends at $q_i$ and winds around $p_i$ counter clock-wise once.

\begin{figure}[h!]
\centering
\def\svgwidth{11cm}
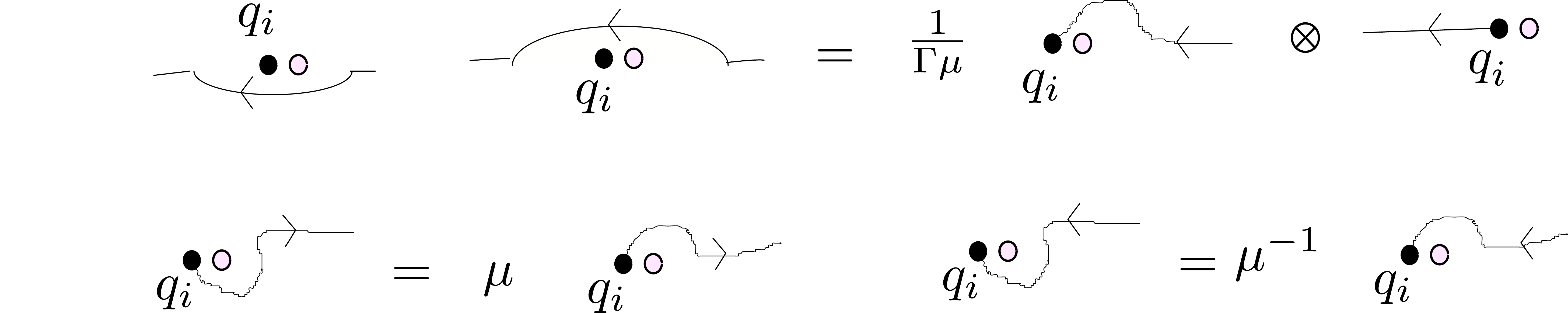
\caption{Derived skein relations}
 \label{deriveskein}
\end{figure}

It should be noted that the relations shown in Figure \ref{deriveskein} can be derived from the ones in Figure \ref{skein12}. And the second relation in Figure \ref{deriveskein} is equivalent to the property that if $\gamma, \gamma' \in \Q_n$ such that $\gamma(0) = q_i$ and $\gamma'(1) = q_i$, then $\gamma_i * \gamma = \mu \gamma, \gamma' * \bar{\gamma_i} = \mu^{-1} \gamma'$, where $*$ means connecting the two adjacent curves, and $\bar{\gamma_i}$ is the curve $\gamma_i$ with reversed direction.

It's not hard to check that any curve can be decomposed into a $($non-commutative$)$ polynomial in $\gamma_{ij}^x \; '$s by repeated applications of the \lq\lq skein" relations. Therefore $\tilde{\A}_{n}$ is generated by $\gamma_{ij}^x \; '$s. Actually they turn out to be free generators after we construct an isomorphism between $\tilde{\A}_{n}$ and $\A_n$ below. Of-course, since $\gamma_{ii}^0 = (1+\mu)\Var$, this doesn't count as part of the free generators.

\begin{figure}[h!]
\centering
\def\svgwidth{5cm}
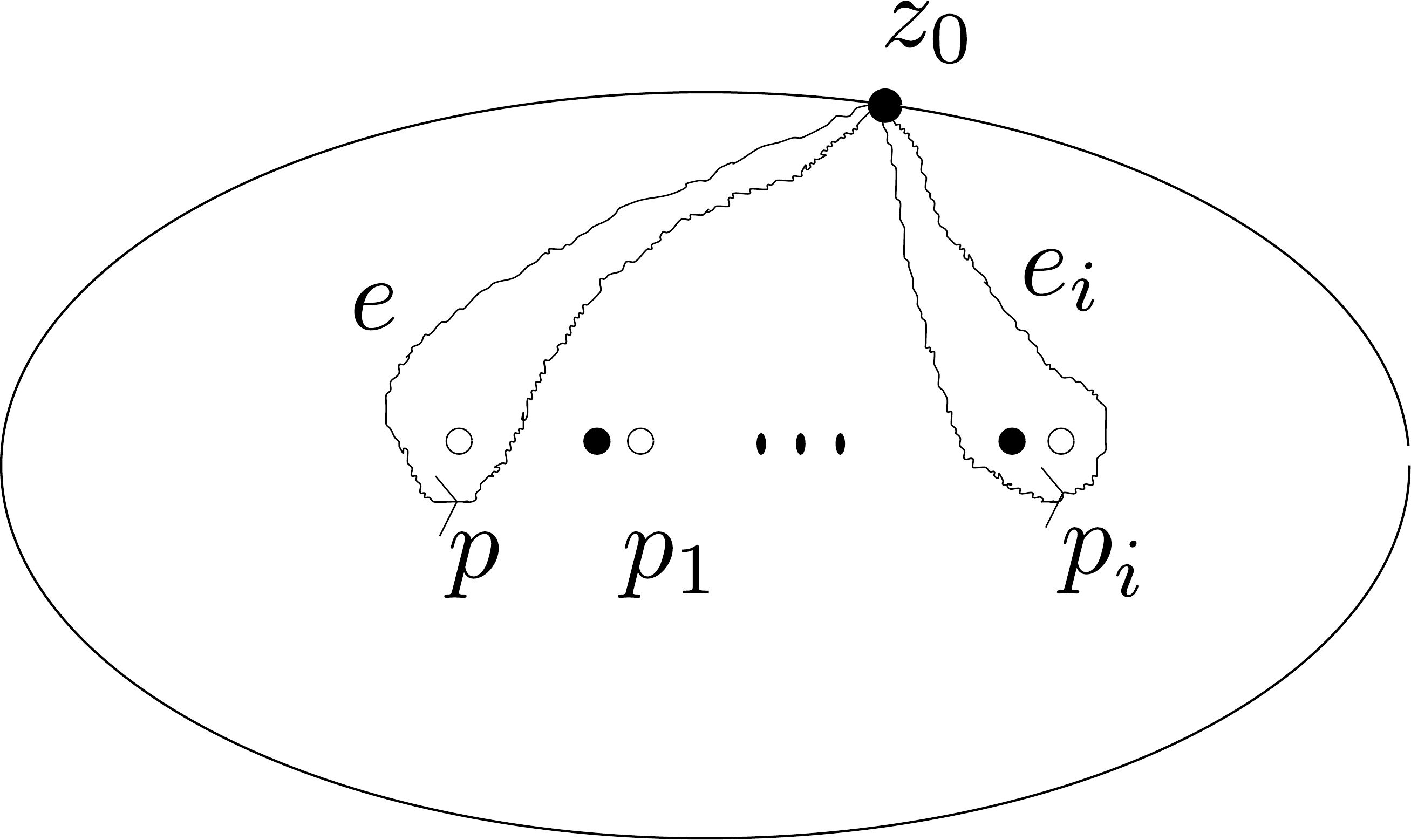
\caption{e and $e_i$ }
\label{ee_i}
\end{figure}

Now pick a base point on the boundary of the disk $D$. To make it explicit, let us pick some $z_0$ on the upper half of the boundary as the base point. The fundamental group of $D_n$ is the free group $F_{n+1} $ on $n+1$ generators, which we denote by $e, e_{1}, \cdots, e_{n}$, where $e_i$ is the loop that winds around $p_i$ counter clock-wise once and $e$ is the loop that winds around $p$ counter clock-wise once. See Figure \ref{ee_i}.

Firstly, we define an intermediate non-commutative algebra $\\$ $\B = R\langle e^{\pm 1}, y_1, y_2, \cdots, y_n \rangle / \mathcal{I}, $ where $\mathcal{I}$ is the two-sided idea generated by $e e^{-1} - 1$, $e^{-1} e - 1$ and $y_i^2 - \Var(1+\mu)y_i, 1 \leq i \leq n$. Define a multiplicative map from $F_{n+1}$ to $\B$ as follows.

$$\tau : F_{n+1} \longrightarrow \B$$

\begin{equation}
\tau(w) =
\begin{cases}
\frac{1}{\Var} y_i - 1    & w = e_i ,1 \leq i \leq n \\
\frac{1}{\Var\mu} y_i - 1    & w = e_i^{-1} ,1 \leq i \leq n \\
e^{\pm 1}              & w = e^{\pm 1}     \\
1                      & w = 1  \\
\end{cases}
\end{equation}

Clearly $\tau(e_i)\tau(e_i^{-1}) = 1 = \tau(1)$ in $\B$. Therefore, we can extend the action of $\tau$ uniquely to arbitrary words to get a well-defined multiplicative map on $F_{n+1}$. Actually $\tau$ extends to an algebra morphism from the group ring $R[F_{n+1}]$ to $\B$.

Next, for $1 \leq i,j \leq n$, we define an $R$-linear map $\alpha_{ij}: R\langle e^{\pm 1}, y_1, y_2, \cdots, y_n \rangle \\ \longrightarrow \A_n$,

$$\alpha_{ij}(e^{i_1}y_{j_1}e^{i_2}y_{j_2} \cdots e^{i_k}y_{j_k}e^{i_{k+1}}) := $$
$$ a_{i,j_1}^{i_1}a_{j_1,j_2}^{i_2} \cdots a_{j_{k-1},j_{k}}^{i_k} a_{j_k,j}^{i_{k+1}}$$

It's easy to check that $\alpha_{ij}$ factors through $\mathcal{I}$ because of the fact that $a_{ii}^0 = (1+\mu)\Var$. Therefore we get an induced map from $\B$ to $\A_n$, which is still denoted by $\alpha_{ij}.$

Finally we can describe the isomorphism between $\tilde{\A}_{n}$ and $\A_n$.

Let $\delta_i$ be the straight line from $z_0$ to $q_i,$ and $\bar{\delta_i}$ be the same line but with reversed direction. For any curve $\gamma \in \Q_n$ with $\gamma(0) = q_i, \gamma(1) = q_j$, let $\tilde{\gamma} = \delta_i*\gamma*\bar{\delta_j}$, then $\tilde{\gamma}$ becomes an element in $\pi_1(D_n,z_0) = F_{n+1}$. Define the isomorphism $\psi : \tilde{\A}_{n} \longrightarrow \A_n$ by $\psi(\gamma):= \alpha_{ij}\tau(\tilde{\gamma})$.

\begin{thm} \label{Atildeiso}
The map $\psi$ defined above is an algebra isomorphism from $\tilde{\A}_{n}$ to $\A_n$ sending $\gamma_{ij}^x$ to $a_{ij}^x$.
\begin{proof}

Clearly, $\psi(\gamma)$ is independent of the choice of $\gamma$ in the equivalence class.

We first show $\psi$ factors through the \lq \lq skein" relations.

It's easy to see that $\psi(\gamma_{ij}^x) = a_{ij}^x$. In particular, $\psi(\gamma_{ii}^0) = a_{ii}^0 = (1+\mu)\Var$, so the first \lq \lq skein" relation is passed to an identity under $\psi$.

Let $C_1, C_2$ denote the two curves passing above and below $p_k$, respectively, in the definition of the second \lq\lq skein" relation in Figure \ref{skein12}. They have the same initial and end points, say $q_i, q_j$. Let $C_3, C_4$ be the curves which ends at $q_k$ and starts at $q_k,$ respectively. So $C_3$ starts from $q_i$ and $C_4$ ends at $q_j$. Let $w_3, w_4$ be the words in $F_{n+1}$ which represent $\tilde{C_3}, \tilde{C_4}$, then it's clear that the words which represent $\tilde{C_1}, \tilde{C_2}$ are $w_3 w_4, w_3 e_k w_4$.

Therefore, $\psi(C_1) + \psi(C_2) = \alpha_{ij}(\tau(w_3)\tau(w_4)) + \alpha_{ij}(\tau(w_3)(\frac{1}{\Var}y_k-1)\tau(w_4)) = \frac{1}{\Var} \alpha_{ij}(\tau(w_3)y_k \tau(w_4)) = \frac{1}{\Var}\alpha_{ik}(\tau(w_3))\alpha_{kj}(\tau(w_4)) = \frac{1}{\Var}\psi(C_3)\psi(C_4)$, which says $\psi$ factors through the second \lq \lq skein" relation.

The above arguments show that $\psi$ is a well-defined algebra morphism. Define the inverse map $\psi': \A_n \longrightarrow \tilde{\A}_{n}$ by sending each $a_{ij}^x$ to $\gamma_{ij}^x$. Noting that $\gamma_{ij}^x$ are generators of $\tilde{\A}_{n}$, it's obvious that $\psi \psi' = Id$ and $\psi' \psi = Id$. Therefore, $\psi$ is an algebra isomorphism.

\end{proof}
\end{thm}

Now we describe a natural action of $\C_n$ on $\tilde{\A}_{n}$.

Recall that the group of isotopy classes of homeomorphisms of $D_n$ with boundary fixed point-wise is the classical braid group on $n+1$ strands $\B_{n+1}.$\footnote{Note that here $D_n$ has $n+1$ punctures.} Here we assume the generators are $\sigma_0, \sigma_1, \cdots, \sigma_{n-1},$ where $\sigma_0$ is the Dehn twist that switches $p$ with $p_1$ counter clock-wise and $\sigma_i$ switches $p_i$ with $p_{i+1}, 1 \leq i \leq n-1.$ Also recall that we identified $\C_n$ with the subgroup of $\B_{n+1}$ which consists of the braids that fix the first puncture. See Section \ref{subsec:markov} for the explicit embedding. Therefore, the elements of $\C_n$ fix the puncture $p$ and permute $\{p_i, 1 \leq i \leq n\}.$ We can furthermore stipulate that the horizontal line segments $p_i q_i$ remain horizontal and of fixed length during the isotopy, so that the elements of $\C_n$ also permute the $q_i\, '$s. It follows that the elements of $\C_n$ act on $\Q_n$. It's also easy to see that this action actually preserves the \lq\lq skein" relations. Therefore, we get a natural action $\tilde{\Phi}$ of $\C_n$ on $\tilde{\A}_{n}$.

\begin{thm}
The algebra isomorphism $\psi : \tilde{\A}_{n} \longrightarrow \A_n$ preserves the action of $\C_n,$ i.e. $\psi\tilde{\Phi}_{\beta} = \Phi_{\beta}\psi,$ for any $\beta \in \C_n.$
\begin{proof}
It suffices to check for any $\beta = \alpha_k$, $\psi\tilde{\Phi}_{\beta} = \Phi_{\beta}\psi$ holds on the generators $\gamma_{ij}^x$.  We left this as an exercise.
\end{proof}
\end{thm}

\begin{rem}
It's worth pointing out that when we want to find the image of some complicated curve in $\tilde{\A}_{n}$ under $\psi,$ it's usually more efficient to use the \lq\lq skein" relations than using the definition directly. Also, instead of memorizing the action of $\C_n$ on the $a_{ij}^x \; '$s, it's much easier to manipulate the \lq\lq skein" relations and the Dehn twists. This provides us another way to calculate the action of a braid $\beta$ on $a_{ij}^x,$ namely, first use a sequence of Dehn twists representing $\beta$ to map $\gamma_{ij}^x$ to some curve, and then decompose this curve into a polynomial of generators using \lq\lq skein" relations, finally replace the generators in the polynomial by the corresponding $a_{ij}^x \; '$s.

For example, to obtain $\Phi_{\alpha_{1}^2 \alpha_0} (a_{12}^{0}),$ we first compute $\tilde{\Phi}_{\alpha_{1}^2 \alpha_0} (\gamma_{12}^{0})$ using Dehn twists that represent $\alpha_{1}^2 \alpha_0.$ See Figure \ref{alpha110}. Then we decompose the resulting curve using \lq\lq skein" relations to get the expression.

$\Phi_{\alpha_{1}^2 \alpha_0} (\gamma_{12}^{0}) = \gamma_{12}^{-1} - \frac{1}{\Var}\gamma_{11}^{-1}\gamma_{12}^{0} + \frac{1}{\Var}\gamma_{12}^{0}\gamma_{22}^{-1} - \frac{1}{\Var^2\mu}\gamma_{12}^{0}\gamma_{21}^{0}\gamma_{12}^{-1} - \frac{1}{\Var^2}\gamma_{12}^{0}\gamma_{21}^{-1}\gamma_{12}^{0}  \\ +\frac{1}{\Var^3\mu}\gamma_{12}^{0}\gamma_{21}^{0}\gamma_{11}^{-1}\gamma_{12}^{0} $

Replacing the $\gamma_{ij}^x\;'$s above with $a_{ij}^x \;'$s, we obtain the expression for $\Phi_{\alpha_{1}^2 \alpha_0} (a_{12}^{0})$.

\begin{figure}[h!]
\def\svgwidth{10cm}
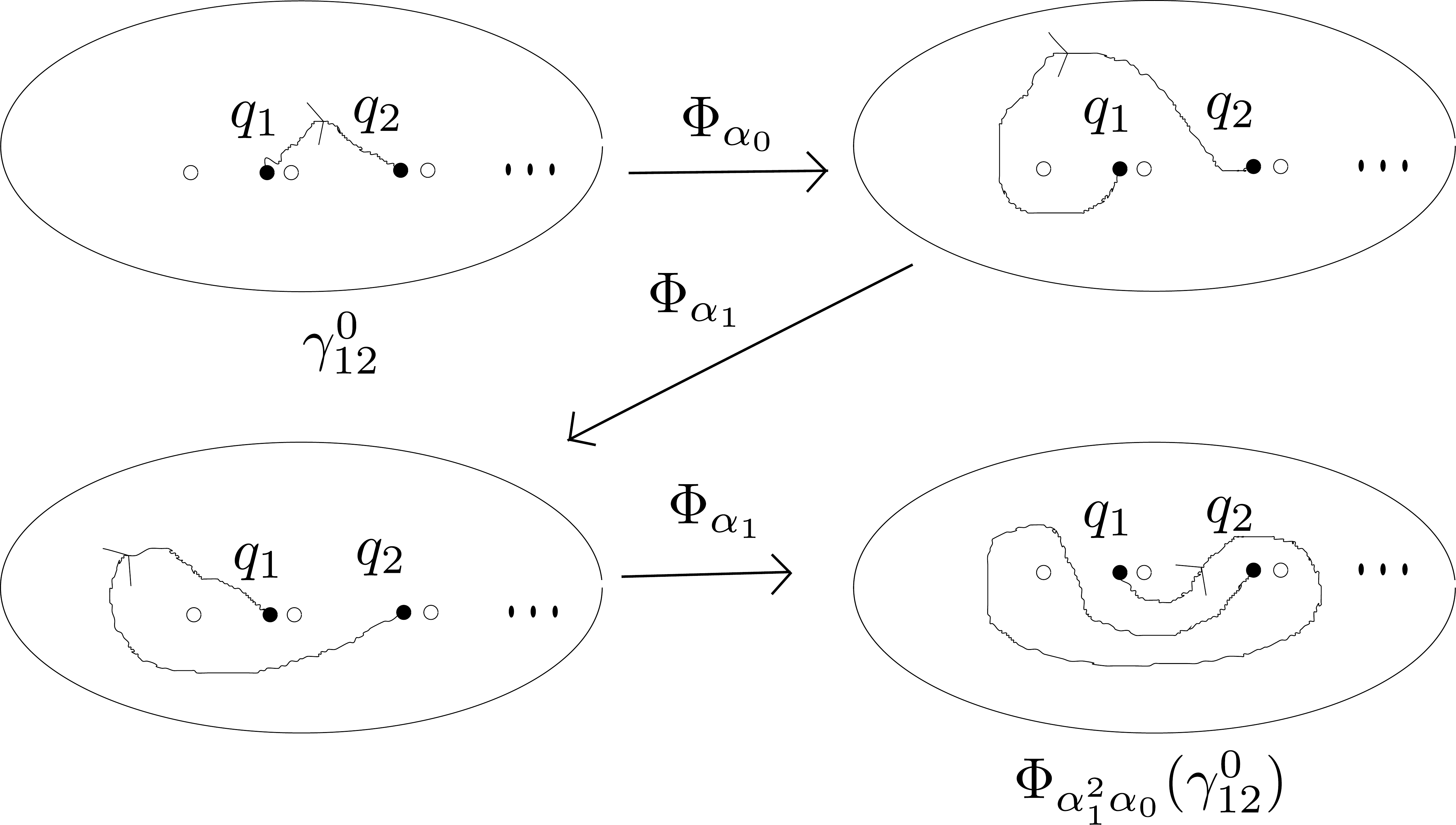
\caption{$\Phi_{\alpha_1^2 \alpha_0} (\gamma_{12}^0)$}
\label{alpha110}
\end{figure}

\end{rem}

There are analogous topological interpretations of the extended actions of $\C_n$ on $\Ap$ and $\Am$.

The procedure goes the same as above, and  we will only point out what modifications should be made at each step.

First of all, let $D_n^{+}$ be the the punctured disk with punctures $p, p_0,p_1, \\ \cdots, p_n$ arranged from left to right and similarly let $D_n^{-}$ be the punctured disk with punctures $p, p_1, \cdots, p_n, p_{n+1}$. Also in both cases, still choose the points $q_i = p_i - \epsilon,$ for some tiny $\epsilon > 0.$ Let $\Q_n^{\pm}$ be the set of equivalence classes of curves in $D_n^{\pm}$ which start and end at the $q_i\; '$s . Define $\tilde{\A}_{n}^{\pm}$ to be the $R$-algebra generated by elements of $\Q_n^{\pm}$ modulo the \lq\lq skein" relations in Figure \ref{skein123}:

\begin{figure}[h!]
\def\svgwidth{11cm}
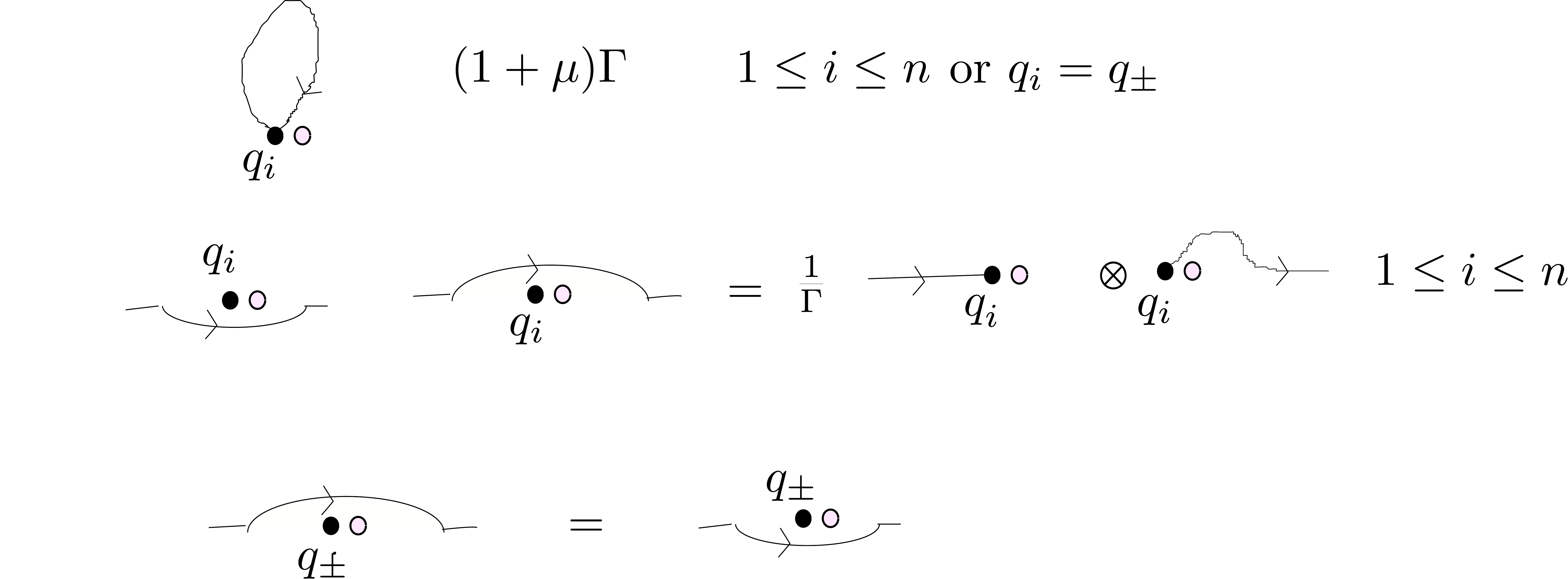
\caption{\textrm{skein relation}}
\label{skein123}
\end{figure}

where $q_{\pm} = q_0$ in the \lq\lq $+$" case and $q_{\pm} = q_{n+1}$ otherwise.

So we added one more relation when defining $\tilde{\A}_{n}^{\pm},$ namely, the curves are allowed to pass through the new puncture $p_0(p_{n+1})$.

The fundamental group of $D_n^{\pm}$ is the free group $F_{n+2}$ generated by $e, e', e_i, 1 \leq i \leq n,$ where $e'$ is the generator that correspond to the new puncture $p_0$ or $p_{n+1}$. We will use the same intermediate algebra $\B$, and the map $\tau$ is extended to $F_{n+2}$ by furthermore defining $\tau(e') = 1.$

In the same way as we defined the isomorphism $\psi$ from $\tilde{\A}_{n}$ to $\A_n$, we can define an isomorphism $\psi^{\pm}$ from $\tilde{\A}_n^{\pm}$ to $\A_n^{\pm}$ which sends $\gamma_{ij}^{x}$ to $a_{ij}^x$.

Next, we extend the action of $\C_n$ to $\tilde{\A}_{n}^{\pm}$.

\begin{figure}[h!]
\def\svgwidth{10cm}
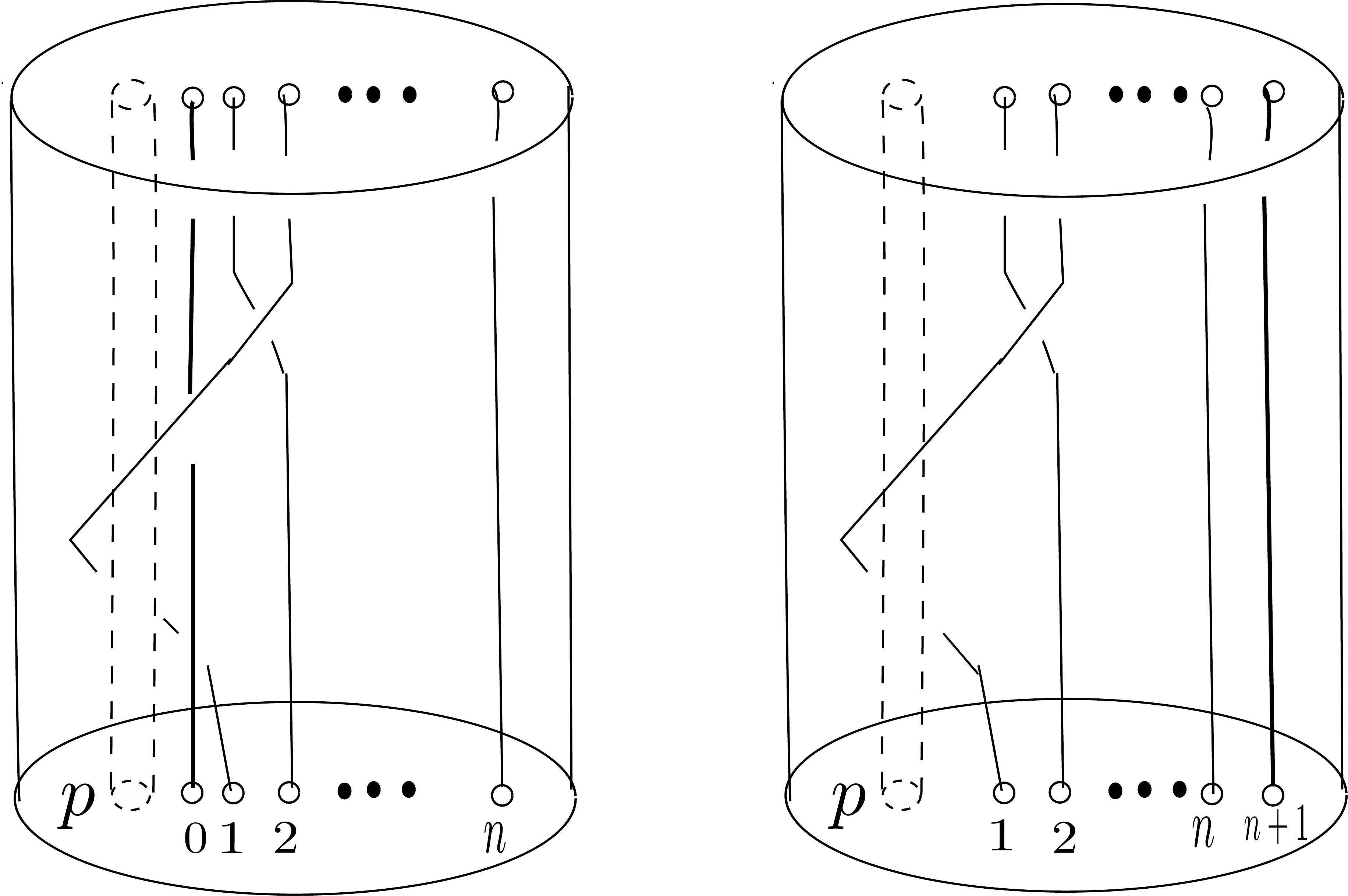
\caption{$\epsilon^{+}(\alpha_1 \alpha_0)$ and $\epsilon^{-}(\alpha_1 \alpha_0)$ }
\label{epsilonpmaction}
\end{figure}

Recall the embedding $\Ep:\C_n \longrightarrow \C_{n+1}$ introduced in Section \ref{subsec:markov}. For notational convenience, we denote the generators of $\C_{n+1}$ by $\alpha_{-1},\alpha_0, \cdots, \alpha_{n-1}.$ Thus the embedding $\Ep$ sends $\alpha_{0}$ to $\alpha_{0}\alpha_{-1}\alpha_{0}$ and $\alpha_i$ to $\alpha_i, \; 1 \leq i \leq n-1$. From the geometrical point of view, $\Ep$ simply inserts a strand labeled by $p_0$ right next to $\{p\} \times [0,1]$.  See the first picture in Figure \ref{epsilonpmaction}.

It's easy to see any braid in $\Ep(\C_n)$ fixes the first two punctures $($the punctures that are labeled by $p$ and $p_0$ $)$. Thus it should be clear that via the embedding $\Ep$, the action of $\C_n$ preserves all the \lq\lq skein" relations defining $\tilde{\A}_{n}^{+}$, and therefore induces an action $\tilde{\Phi}^{+}$ on $\tilde{\A}_{n}^{+}$.

For the action $\tilde{\Phi}^{-}$ of $\C_n$ on $\tilde{\A}_{n}^{-}$, we use the other embedding $\Em: \C_n \longrightarrow \C_{n+1}$. Note that here the generators of $\C_{n+1}$ are $\alpha_0, \cdots , \alpha_{n},$ and $\Em(\alpha_i) = \alpha_i, 0 \leq i \leq n-1.$ And the map $\Em$ inserts a strand labeled by $p_{n+1}$ on the right of the braid. See the second picture in Figure \ref{epsilonpmaction}.

Here in Figure \ref{epsilonpmaction} we use $i$ to represent $p_i$.

Again, since elements of $\Em(\C_n)$ fix $p_{n+1}$, they preserve the  \lq\lq skein" relations that define $\tilde{\A}_n^{-}$. We thus get an induced action $\tilde{\Phi}^{-}$ of $\C_n$ on $\tilde{\A}_n^{-}$.

$\tilde{\A}_{n}$ can obviously be embedded as a subalgebra into $\tilde{\A}_{n}^{\pm}$. We have the following theorem which relates the topological interpretations of the actions of $\C_n$ to the algebraic interpretations.

\begin{thm}\label{Atildepmiso}
The maps $\psi^{\pm}: \tilde{\A}_{n}^{\pm} \longrightarrow \A_n^{\pm}$ are algebra isomorphisms and commute with the extended actions of $\C_n,$ namely, for any $\beta \in \C_n$, $\psi^{\pm}\tilde{\Phi}^{\pm}_{\beta} = \Phi^{\pm}_{\beta}\psi^{\pm}.$ Moreover, $(\tilde{\Phi}^{\pm}_{\beta})_{| \tilde{\A}_{n}} = \tilde{\Phi}_{\beta}$, and the following diagram commutes:

\begin{equation}
\xymatrix{
\tilde{\A}_{n} \ar[r]^-{\psi} \ar@{^{(}->}[d]  &  \A_n \ar@{^{(}->}[d] \\
\tilde{\A}_{n}^{\pm} \ar[r]^-{\psi^{\pm}}          &  \A_n^{\pm} \\
}
\end{equation}

And each of the maps in the above diagram preserves the action of $\C_n$.
\begin{proof}
Proofs are analogous to that of Theorem \ref{Atildeiso}.
\end{proof}
\end{thm}

It is worth noting that the actions of $\tilde{\Phi}^{\pm} $ on $\tilde{\A}_{n}^{\pm}$ and the action of $\tilde{\Phi}$ on $\tilde{\A}_{n}$ can also be visualized as follows.

For a braid $\beta \in \C_n,$ draw a braid diagram of $\beta$ inside $D_n \times [0,1]$, such that the intersection of the braid with $D_n \times \{0,1\}$ are exactly the punctures $p_i\;'$s. Perturb the braid diagram to get a parallel copy of it such that the intersection of the copy with $D_n \times \{0,1\}$ are the $q_i\;'$s. For any curve $\gamma \subset D_n \times \{0\}$ representing some element in $\tilde{\A}_n$, slide $\gamma$ along the copy diagram in the complement of the braid diagram until it reaches $D_n \times \{1\}$, then the resulting curve is $\tilde{\Phi}_{\beta}(\gamma)$.

To visualize $\tilde{\Phi}^{\pm}_{\beta}$, we draw a braid diagram of $\Epm(\beta)$ inside $D_n^{\pm} \times I$, make a parallel copy of it, and slide any curve along the copy diagram up to $D_n^{\pm} \times \{1\}.$

With the above observations, we have the following simple but important proposition.
\begin{prop} \label{prop:curve connect}
 Let $\beta \in \C_n$ be a braid.

$1)$. For any two curves $\gamma_1,\gamma_2 \in \Q_n^{+},$ such that $\gamma_1(1) = \gamma_2(0) = q_0$ and $\gamma_1(0) = q_i, \gamma_2(1) = q_j$ for some $1 \leq i,j \leq n$, then $\gamma_1*\gamma_2$ is a curve in $\Q_n$ from $q_i$ to $q_j$, and we have $\tilde{\Phi}_{\beta}(\gamma_1*\gamma_2) = \tilde{\Phi}^{+}_{\beta}(\gamma_1)*\tilde{\Phi}^{+}_{\beta}(\gamma_2)$.

$2)$. For any two curves $\gamma_1,\gamma_2 \in \Q_n^{-},$ such that $\gamma_1(1) = \gamma_2(0) = q_{n+1}$ and $\gamma_1(0) = q_i, \gamma_2(1) = q_j$ for some $1 \leq i,j \leq n$, then $\gamma_1*\gamma_2$ is a curve in $\Q_n$ from $q_i$ to $q_j$, and we have $\tilde{\Phi}_{\beta}(\gamma_1*\gamma_2) = \tilde{\Phi}^{-}_{\beta}(\gamma_1)*\tilde{\Phi}^{-}_{\beta}(\gamma_2)$.

Here, and throughout the paper, $*$ again means connecting the two curves.
\end{prop}

\begin{rem} \label{rem picture action}
From now on, we will identify $\tilde{\A}_{n}$ with $\A_n$, $\tilde{\A}_{n}^{\pm} $ with $\Apm$, $\gamma_{ij}^x $ with $a_{ij}^x$ via the corresponding isomorphisms and identify $\tilde{\Phi}_{\beta}$ with $\Phi_{\beta}$, $\tilde{\Phi}^{\pm}_{\beta}$ with $\Phipm{\beta},$ respectively. A useful picture to keep in mind is as follows. $a_{ij}^x$ is the left arc diagram described in Figure \ref{gammaij}. The action $\Phi_{\beta}\, ( \Phipm{\beta})$ of $\beta$ on some curve is to slide that curve along the parallel copy of the braid diagram that represents $\beta \, (\Epm(\beta))$ up to $D_n \times \{1\} \, (D_n^{\pm} \times \{1\}).$
\end{rem}

We can also define the \lq\lq $*$" operation on some elements of $\A_n$.

\begin{definition} \label{def:connect}

$1).$ Let $P,Q \in \Ap$ such that $P= \sum\limits_{x \in \Z}\sum\limits_{i=1}^n P_i^x a_{i0}^x, Q= \sum\limits_{y \in \Z}\sum\limits_{j=1}^n a_{0j}^y Q_j^y, \, P_i^x, Q_j^y \in \A_n$, then $P*Q \in \A_n$ is defined to be $\sum\limits_{x,y \in \Z}\sum\limits_{i,j=1}^n P_i^x a_{ij}^{x+y}Q_j^y.$

$2).$ Let $P,Q \in \Am$ such that $P= \sum\limits_{x \in \Z}\sum\limits_{i=1}^n P_i^x a_{i,n+1}^x, Q= \sum\limits_{y \in \Z}\sum\limits_{j=1}^n a_{n+1,j}^y Q_j^y, \, \\ P_i^x, Q_j^y \in \A_n$, then $P*Q \in \A_n$ is defined to be $\sum\limits_{x,y \in \Z}\sum\limits_{i,j=1}^n P_i^x a_{ij}^{x+y}Q_j^y.$

$3).$ Two elements $P, Q \in A_{n}^{\pm}$ are called connectable, if they satisfy the condition in one of the above two definitions.
\end{definition}

\begin{prop} \label{connectable}
If $P,Q \in A_{n}^{\pm}$ are connectable, then for any $\beta \in \C_n$, $\Phipm{\beta}(P), \Phipm{\beta}(Q)$ are also connectable, and $\Phi_{\beta}(P*Q) = \Phipm{\beta}(P)*\Phipm{\beta}(Q)$.
\begin{proof}
We will only prove the \lq\lq +" case. The proof of the other case is analogous.

Let $P,Q$ be as described in $1)$ of Definition \ref{def:connect}, then for $\beta \in \C_n$, $\Phip{\beta}(P) = \sum\limits_{x \in \Z}\sum\limits_{i=1}^n \Phip{\beta}(P_i^x) \Phip{\beta}(a_{i0}^x) = \sum\limits_{x \in \Z}\sum\limits_{i=1}^n \Phi_{\beta}(P_i^x) \Phip{\beta}(a_{i0}^x),$ and similarly, $\Phip{\beta}(Q) = \sum\limits_{y \in \Z}\sum\limits_{j=1}^n \Phip{\beta}(a_{0j}^y) \Phi_{\beta}(Q_j^y) $. Clearly, $\Phip{\beta}(a_{i0}^x)$ and $\Phip{\beta}(a_{0j}^y)$ are connectable, so $\Phip{\beta}(P)$ and $\Phip{\beta}(Q)$ are connectable. Moreover,

$\Phip{\beta}(P) * \Phip{\beta}(Q) = \sum\limits_{x,y \in \Z}\sum\limits_{i,j=1}^n \Phi_{\beta}(P_i^x) \{\Phip{\beta}(a_{i0}^x)*\Phip{\beta}(a_{0j}^y)\}\Phi_{\beta}(Q_j^y) \\ \xlongequal{Proposition \, \ref{prop:curve connect}} \sum\limits_{x,y \in \Z}\sum\limits_{i,j=1}^n \Phi_{\beta}(P_i^x) \Phi_{\beta}(a_{ij}^{x+y})\Phi_{\beta}(Q_j^y) = \\ \Phi_{\beta}(\sum\limits_{x,y \in \Z}\sum\limits_{i,j=1}^n P_i^x a_{ij}^{x+y}Q_j^y) = \Phi_{\beta}(P*Q). $

\end{proof}
\end{prop}

\section{The framed knot invariant}
\label{sec:knot invariant}
From now on, we will assume the closure of $\beta \in C_n$ is a knot in $\Lens$.

In this section, first we give the definition of the framed knot invariant. Since the knot invariant looks very complicated at first glance, we will compute some examples after the definition. We then proceed to give some ancillary results, and finally prove the invariance under Markov moves.

\subsection{Definition of the invariant}
\label{subsec:def invariant}

Here are some notations we will use to define the invariant.

Let $M_{\infty}(\A_n)$ denote the set of $\infty \times \infty$ matrices with elements in $\A_n$, namely, the rows and columns of a matrix in $M_{\infty}(\A_n)$ are both indexed by integers. We call a matrix {\it{row-finite}} if there are only finitely many non-zero entries in each row. A {\it{column-finite}} matrix is defined analogously. If $M,N $ are two matrices in $\M_{\infty}(\A_n)$, in general the multiplication of them is not well-defined. However, if $M$ is {\it{row-finite}} or $N$ is {\it{column-finite}}, then $M\cdot N$ is well-defined. And the associativity is satisfied whenever multiplications make sense. All throughout the paper, the matrices always satisfy the above condition when they are multiplied together, and for $x,y \in \Z$, we will use $M^{xy}$ to refer to the $(x,y)$-entry of $M$. We will also use an element $c \in A_n$ to represent the scalar matrix in $M_{\infty}(\A_n)$ which has entry $c$ on the diagonal and $0$ elsewhere. Let $M_n(M_{\infty}(\A_n))$ denote the set of $n \times n$ matrices with entries in $M_{\infty}(\A_n)$.

Recall $\Epm: \C_n \longrightarrow \C_{n+1}$ are the two embeddings, and for $\beta \in \C_n$, $(\Phim{\beta})_{|\A_n} = \Phi_{\beta} = (\Phip{\beta})_{|\A_n}.$

It's not hard see $($perhaps easier from the topological interpretation$)$ that for $1 \leq i \leq n, x \in \Z,$ $\Phim{\beta}(a_{i,n+1}^x)$ can be written as a finite linear combinations of $a_{k,n+1}^z, 1 \leq k \leq n,  z \in \Z $ with coefficients in $\A_n$. A similar argument holds for $\Phim{\beta}(a_{n+1,i}^x), \Phip{\beta}(a_{i,0}^x), \Phip{\beta}(a_{0,i}^x)$. For example, $\Phip{\beta}(a_{0,i}^x)$ is a finite linear combinations of $a_{0,k}^z$ with coefficients in $\A_n$ multiplied on the right. Explicitly, this is how we define $\Phiml{\beta}, \Phimr{\beta}, \Phipl{\beta} , \Phipr{\beta} \in M_n(M_{\infty}(\A_n))$.

For each $\beta \in \C_n, 1 \leq i,j \leq n, x, y \in \Z$, define

$$\Phim{\beta}(a_{i,n+1}^{x}) = \sum\limits_{k=1}^{n} \sum\limits_{z \in \Z} (\Phiml{\beta})_{ik}^{xz}a_{k,n+1}^z$$
$$\Phim{\beta}(a_{n+1,j}^{y}) = \sum\limits_{k=1}^{n} \sum\limits_{z \in \Z} a_{n+1,k}^z(\Phimr{\beta})_{kj}^{zy}$$
$$\Phip{\beta}(a_{i,0}^{x}) = \sum\limits_{k=1}^{n} \sum\limits_{z \in \Z} (\Phipl{\beta})_{ik}^{xz}a_{k,0}^z$$
$$\Phip{\beta}(a_{0,j}^{y}) = \sum\limits_{k=1}^{n} \sum\limits_{z \in \Z} a_{0,k}^z(\Phipr{\beta})_{kj}^{zy}$$

where $(\Phiml{\beta})_{ik}^{xz}$ is the $(x,z)$-entry of the $\infty \times \infty$ matrix $(\Phiml{\beta})_{ik}$ which is the $(i,k)$-entry of the $n \times n$ matrix $\Phiml{\beta}$. So we have $\Phiml{\beta} \in M_n(M_{\infty}(\A_n)).$ A similar statement holds for the other three symbols.

Define $a_{ij} \in M_{\infty}(\A_n)$ by $(a_{ij})^{xy} = a_{ij}^{x+y}$ and define $A \in M_n(M_{\infty}(\A_n))$ by $A_{ij} = a_{ij}.$

\begin{lem}
For $\beta \in \C_n, 1 \leq i,j \leq n,$ $(\Phiml{\beta})_{ij},(\Phipl{\beta})_{ij}$ are {\it{row-finite}} and $(\Phimr{\beta})_{ij},(\Phipr{\beta})_{ij}$ are {\it{column-finite}}.
\begin{proof}
These are direct consequences of the definitions.
\end{proof}
\end{lem}

\begin{rem}
Actually, $(\Phipl{\beta})_{ij},(\Phipr{\beta})_{ij}$ are both {\it{row-finite}} and {\it{column-finite}}. This is due to a careful inspection of the action $\Phip{\beta}.$ We will not use this property though.
\end{rem}

For $1 \leq p,q \leq n, f \in \Z$, let $\Lambda_{f;p,q} \in M_n(M_{\infty}(\A_n))$ be the diagonal matrix with the $(p,p)$-th entry $\lambda$, the $(q,q)$-entry $\mu^{-f}$ and other diagonal entries $1$.

\begin{definition} \label{def invariant}
Let $\beta \in \C_n, 1 \leq p, q \leq n, f \in \Z$, then $HC_0(\beta;f; p, q)$ is defined to be the $R$-algebra $\A_n$ modulo the two sided idea $\mathcal{I}_{\beta;f;p,q}$ generated by the entries of the entries of following matrices:

$$A - \Lambda_{f;p,q}\Phiml{\beta}A$$
$$A - A \Phimr{\beta}\Lambda_{f;p,q}^{-1}$$
$$A - \Lambda_{f;p,q}\Phipl{\beta}A$$
$$A - A \Phipr{\beta}\Lambda_{f;p,q}^{-1}$$
\end{definition}

\begin{rem}
\begin{enumerate}
\item For a matrix $M \in M_n(M_{\infty}(\A_n))$, the phrase \lq\lq the entries of the entries of M" is really awkward. From now on, we will use \lq\lq the elements of $M$" to stand for \lq\lq the entries of the entries of $M$".

\item  Note that $(\Lambda_{p,q;f}\Phiml{\beta}A)_{ij}^{xy} = \sum\limits_{k=1}^{n}\sum\limits_{z \in \Z}\lambda^{\delta_{i,p}}\mu^{-f \delta_{i,q}}(\Phiml{\beta})_{ik}^{xz}A_{kj}^{zy} \\
     = \sum\limits_{k=1}^{n}\sum\limits_{z \in \Z}\lambda^{\delta_{i,p}}\mu^{-f \delta_{i,q}}((\Phiml{\beta})_{ik}^{xz}a_{k,n+1}^z) * a_{n+1,j}^y \\ = \lambda^{\delta_{i,p}}\mu^{-f \delta_{i,q}}\Phim{\beta}(a_{i,n+1}^x) * a_{n+1,j}^y$.

    Since $A_{ij}^{xy} = a_{ij}^{x+y} = a_{i,n+1}^x*a_{n+1,j}^{y},$ the relations in $\mathcal{I}_{\beta;f;p,q}$ are the same as the following:

    $a_{i,n+1}^x*a_{n+1,j}^y - \lambda^{\delta_{i,p}}\mu^{-f \delta_{i,q}}\Phim{\beta}(a_{i,n+1}^x)*a_{n+1,j}^y$,

    $a_{i,n+1}^x*a_{n+1,j}^y - \lambda^{-\delta_{j,p}}\mu^{f \delta_{j,q}}a_{i,n+1}^x*\Phim{\beta}(a_{n+1,j}^y)$,

    $a_{i,0}^x*a_{0,j}^y - \lambda^{\delta_{i,p}}\mu^{-f \delta_{i,q}}\Phip{\beta}(a_{i,0}^x)*a_{0,j}^y$,

    $a_{i,0}^x*a_{0,j}^y - \lambda^{-\delta_{j,p}}\mu^{f \delta_{j,q}}a_{i,0}^x*\Phip{\beta}(a_{0,j}^y), \, \forall 1 \leq i,j \leq n,\,  x,y \in \Z$.

\end{enumerate}
\label{rem relation}
\end{rem}

Note that for $\beta \in C_n$, it has a natural action by permutation on the set $\{1, \cdots, n\}$. Our convention here is that the braid diagram always goes upward, and if the $i$-th strand ends at the $j$-th position, then $\beta(i)=j$.

\begin{lem}
For $\beta \in \C_n, 1 \leq p, q \leq n, f \in \Z$, we have $HC_0(\beta;f; p,q) \simeq HC_0(\beta;f; \beta(p),q) \simeq HC_0(\beta;f; p,\beta(q))$.
\begin{proof}
Define $\psi: HC_0(\beta;f; p,q) \longrightarrow HC_0(\beta;f; \beta(p),q)$ by $\psi(a_{ij}^x) = \lambda^{-\delta_{i,\beta(p)}}a_{ij}^x\lambda^{\delta_{j,\beta(p)}}$. We need to check that $\psi$ sends $\mathcal{I}_{\beta;f;p,q}$ to $\mathcal{I}_{\beta;f;\beta(p),q}$.

Note that $\Phim{\beta}(a_{i,n+1}^x)*a_{n+1,j}^y$ can be written as a non-commutative polynomial in which each monomial is of the form $a_{\beta(i),i_1}^{x_1}a_{i_1,i_2}^{x_2}\cdots a_{i_k,j}^{x_{k+1}}$, thus we have

$\psi(a_{ij}^{xy} - \lambda^{\delta_{i,p}}\mu^{-f \delta_{i,q}}\Phim{\beta}(a_{i,n+1}^x)*a_{n+1,j}^y) \\
= \lambda^{-\delta_{i,\beta(p)}}a_{ij}^{xy}\lambda^{\delta_{j,\beta(p)}} - \lambda^{\delta_{i,p}}\mu^{-f \delta_{i,q}}\lambda^{-\delta_{\beta(i),\beta(p)}}\Phim{\beta}(a_{i,n+1}^x)*a_{n+1,j}^y \lambda^{\delta_{j,\beta(p)}}\\
= \lambda^{-\delta_{i,\beta(p)}}(a_{ij}^{xy} - \lambda^{\delta_{i,\beta(p)}}\mu^{-f \delta_{i,q}}\Phim{\beta}(a_{i,n+1}^x)*a_{n+1,j}^y) \lambda^{\delta_{j,\beta(p)}} \in \mathcal{I}_{\beta;f;\beta(p),q}$.

The other three relations in $\mathcal{I}_{\beta;f;p,q}$ can be shown analogously that they are mapped under $\psi$ to $\mathcal{I}_{\beta;f;\beta(p),q}$. Thus the map $\psi$ is well defined. It's also clear that it's an isomorphism.

The isomorphism $HC_0(\beta;f; p,q)  \simeq HC_0(\beta;f; p,\beta(q))$ can be proved similarly by sending $a_{ij}^x $ to $\mu^{k\delta_{i,\beta(p)}}a_{ij}^x\mu^{-k\delta_{j,\beta(p)}}$.
\end{proof}
\end{lem}

\begin{cor}
If the closure of $\beta \in \C_n$ is a knot in $\Lens$, then $HC_0(\beta;f; p,q)$ is independent of the values of $p,q$.
\end{cor}

Apparently from the definition, $HC_0(\beta;f;p,p)$ can be obtained from $HC_0(\beta;0;p,p)$ by replacing $\lambda$ by $\lambda\mu^{-f}$. We will use the notations $HC_0(\beta;f;p) = HC_0(\beta;f;p,p)$, $HC_0(\beta;f) = HC_0(\beta;f;1,1) $ and $HC_0(\beta) = HC_0(\beta;0;1,1) $. By the corollary above, $HC_0(\beta;f;p)$ is dependent of the choice of $p$, so we have $HC_0(\beta;f) \simeq HC_0(\beta;f;p)$ for any $p$.

The following theorem is our main result.

\begin{thm} \label{HC0 braid}
Let $\beta, \alpha \in C_n, f \in \Z$ such that the closure of $\beta$ in $\Lens$ is a knot, then we have the following isomorphisms:

$1).\, HC_0(\beta;f) \simeq HC_0(\alpha^{-1}\beta\alpha;f)$;

$2).\, HC_0(\beta;f) \simeq HC_0(\Em(\beta)\alpha_n;f-1) \simeq HC_0(\Em(\beta)\alpha_n^{-1};f+1)$;

$3).\, HC_0(\beta;f) \simeq HC_0(\Ep(\beta)\alpha_n;f-1) \simeq HC_0(\Ep(\beta)\alpha_n^{-1};f+1)$;
\end{thm}

We will give a proof of the theorem in Section \ref{subsec:invariace proof}.

Endow $\Lens$ with the standard orientation. Let $K$ be a framed oriented knot in $\Lens$ with $l,m$ the homotopy classes of the longitude and the meridian of $K$ in $\pi_1(\Lens \setminus K)$. The orientations of $K$ and $\Lens$ determine the meridian class $m$ uniquely. More precisely, let $\nu(K)$ be the tubular neighborhood of $K$, which is homeomorphic to $K \times D^2$. Choose an orientation on $D^2$ so that the homeomorphism of $\nu(K)$ with $K\times D^2$ is orientation preserving. Then for any $z \in K$, the image of $z \times \partial{D^2}$ under the homeomorphism determines the meridian class. Assume $K$ is represented by the closure of a braid $\beta \in \C_n$, and that $[l] = [\hat{\beta'}][m]^f$, where $\beta'$ is a parallel push-off copy of $\beta$, then $HC_0(K;l)$ is defined to be $HC_0(\beta;f)$.

\begin{cor}
$HC_0(K;l)$ as an $R$-algebra is a framed knot invariant for knots in $\Lens.$
\begin{proof}
For a braid diagram $\beta \in \C_n$, let $\beta'$ be the parallel push-off copy of $\beta$. Then we have $[\hat{\beta'}][m]^{\pm 1} = [\widehat{(\Ep(\beta)\alpha_n^{\pm 1})'}], \, [\hat{\beta'}][m]^{\pm 1} = [\widehat{(\Em(\beta)\alpha_n^{\pm 1})'}]$ and for any $\alpha \in \C_n$, we have $[\hat{\beta'}] = [\widehat{(\alpha^{-1}\beta\alpha)'}].$
\end{proof}
\end{cor}

\begin{rem}
The invariant $HC_0(K;l)$ is conjectured to be the $0$-th knot contact homology of $K$, which is defined to be the $0$-th Legendrian contact homology of $\Lambda_K$ in $ST^{*}(\Lens),$ where $ST^{*}(\Lens)$ is the unit cotangent bundle of $\Lens$ and $\Lambda_K$ is the unit conormal bundle of $K$. As this paper is not relevant to proving this conjecture, readers should just treat $HC_0$ purely as a name.
\end{rem}

\subsection{Examples} \label{subsec:example}

Before proving invariance, we first look at some examples.

\begin{example}

$1).$ \textbf{Unknot}. The most simple example is the unknot represented by the identity element $e$ in $\C_1$. We compute $HC_0(e;f)$ for $f \in \Z$. In this case, it's clear that $\Phipl{e},\Phipr{e},\Phiml{e},\Phimr{e}$ are all identity matrices, thus all the relations in $\mathcal{I}_{e;f;1,1}$ become $(1 - \lambda\mu^{-f})a_{11}^x $, and so $HC_0(e;f) \simeq R\langle a_{11}^x, x \in \Z\rangle / \langle (1 - \lambda\mu^{-f})a_{11}^x \rangle$.

$2).$ \bm{$\widehat{\alpha_0^2}$}. Set $\beta = \alpha_0^2, \Lambda = \Lambda_{\beta;0;1,1}$. We first compute $\Phipl{\beta}, \Phipr{\beta}$. Direct calculations show that $\Phip{\beta}(a_{10}^x) = \mu^2 a_{10}^{x-2}, \, \Phip{\beta}(a_{01}^y) = \mu^{-2}a_{01}^{y+2}.$ Thus we have $(\Phipl{\beta})_{11}^{xy} = \mu^2\delta_{x-2,y}, \, (\Phipr{\beta})_{11}^{xy} = \mu^{-2}\delta_{x-2,y},$ and therefore $(\Lambda\Phipl{\beta}A)_{11}^{xy} = \lambda\mu^2 a_{11}^{x+y-2}, \, (A\Phipr{\beta}\Lambda^{-1})_{11}^{xy} = (\lambda\mu^2)^{-1}a_{11}^{x+y+2}.$ So the the third and fourth relation defining $\mathcal{I}_{\beta;0;1,1}$ both are $a_{11}^{x+2} - \lambda\mu^2 a_{11}^{x}, x \in \Z$.

Now we compute $\Phiml{\beta}, \Phimr{\beta}$. By definition, $\Phim{\alpha_0}(a_{11}^x) = a_{11}^x, \, \Phim{\alpha_0}(a_{12}^x) = -\mu a_{12}^{x-1} + \frac{1}{\Var}a_{11}^{x}a_{12}^{-1}$. Therefore,

$\Phim{\alpha_0^2}(a_{12}^x) = -\mu \Phim{\alpha_0}(a_{12}^{x-1}) + \frac{1}{\Var}\Phim{\alpha_0}(a_{11}^{x})\Phim{\alpha_0}((a_{12}^{-1}) \\ = \mu^{2}a_{12}^{x-2} - \frac{\mu}{\Var}a_{11}^{x-1}a_{12}^{-1} - \frac{\mu}{\Var}a_{11}^xa_{12}^{-2} + \frac{1}{\Var^2}a_{11}^xa_{11}^{-1}a_{12}^{-1}. \\ $

By Part $(2)$ of Remark \ref{rem relation},

$(\Lambda\Phiml{\beta}A)_{11}^{xy}-A_{11}^{xy} = \lambda(\mu^{2}a_{11}^{x+y-2} - \frac{\mu}{\Var}a_{11}^{x-1}a_{11}^{y-1} - \frac{\mu}{\Var}a_{11}^xa_{11}^{y-2} + \frac{1}{\Var^2}a_{11}^x a_{11}^{-1}a_{11}^{y-1}) - a_{11}^{x+y}$.

Similarly,

$(A\Phimr{\beta}\Lambda^{-1})_{11}^{xy}-A_{11}^{xy} = (\lambda\mu^2)^{-1}(a_{11}^{x+y+2} - \frac{1}{\Var}a_{11}^{x+1}a_{11}^{y+1}-\frac{1}{\Var}a_{11}^{x+2}a_{11}^{y} + \frac{1}{\Var^2}a_{11}^{x+1}a_{11}^{1}a_{11}^{y}) - a_{11}^{x+y}$.

Since we have $a_{11}^{x+2} - \lambda\mu^2 a_{11}^{x}$, then the above two relations can be simplified as

$a_{11}^{x-1}a_{11}^{y-1} + a_{11}^xa_{11}^{y-2} - \frac{1}{\Var \mu}a_{11}^x a_{11}^{-1}a_{11}^{y-1}$ and

$a_{11}^{x-1}a_{11}^{y-1} + a_{11}^xa_{11}^{y-2} - \frac{1}{\Var }a_{11}^{x-1} a_{11}^{1}a_{11}^{y-2}$.

And only parities of $x$ and $y$ will make a difference in the above two relations.

Direct calculation shows that  $HC_0(\beta) \simeq R[X]/\langle (1-\mu)X, X^2 - \Var^2\lambda(1+\mu)^2\rangle$. Replacing $\lambda$ by $\lambda\mu^{-f}$, we obtain $HC_0(\beta;f)$.

It will be shown in Section \ref{subsec:torus knot} that $\widehat{\alpha_0^2}$ is a particular knot in a large family of knots, namely the torus knots. Explicitly, it is the $(1,2)$-knot. See Section \ref{subsec:torus knot} for a definition of torus knots and more examples.

\end{example}

\subsection{Properties of $\Phi^{\pm L}, \Phi^{\pm R}$}
\label{subsec:ancilary}

We give several propositions which will be used in proving the invariance of  $HC_0(K;l)$. A similar version of these propositions are proved in \cite{ng2005knotI} where the author defined the $HC_0$ invariant for knots in $S^3$.

If $\phi$ is an algebra morphism from $\A_n$ to $\A_n$, and $M \in M_n(M_{\infty}(\A_n))$,  we denote by $\phi(M)$ or $M(\phi)$ the matrix obtained from $M$ by replacing each $a_{ij}^x$ by $\phi(a_{ij}^x)$. Recall in last subsection, we defined the four matrices $\Phiml{\beta}, \Phimr{\beta}, \Phipl{\beta},\Phipr{\beta} \in M_n(M_{\infty}(\A_n))$ for $\beta \in \C_n$.

\begin{prop}\label{Phi(beta 1 2)}
Let $\beta_1, \beta_2 \in \C_n$ be two braids, then we have
$$\Phiml{\beta_1\beta_2} = \Phiml{\beta_2}(\Phi_{\beta_1} )\Phiml{\beta_1}$$
$$\Phimr{\beta_1\beta_2} = \Phimr{\beta_1}\Phimr{\beta_2}(\Phi_{\beta_1} )$$
$$\Phipl{\beta_1\beta_2} = \Phipl{\beta_2}(\Phi_{\beta_1} )\Phipl{\beta_1}$$
$$\Phipr{\beta_1\beta_2} = \Phipr{\beta_1}\Phipr{\beta_2}(\Phi_{\beta_1} )$$

\begin{proof}
The proof of the four equalities are straight forward and completely analogous, so we will just prove the first one.

By definition, $\Phim{\beta_2}(a_{i,n+1}^x) =  \sum\limits_{k=1}^{n} \sum\limits_{z \in \Z} (\Phiml{\beta_2})_{ik}^{xz}a_{k,n+1}^z$. So

$$\Phim{\beta_1\beta_2}(a_{i,n+1}^x) = \Phim{\beta_1}\Phim{\beta_2}(a_{i,n+1}^x)$$
$$= \sum\limits_{k=1}^{n} \sum\limits_{z \in \Z} \Phim{\beta_1}((\Phiml{\beta_2})_{ik}^{xz})\Phim{\beta_1}(a_{k,n+1}^z) $$
$$= \sum\limits_{k,j=1}^{n} \sum\limits_{z,y \in \Z} \Phiml{\beta_2}(\Phi_{\beta_1} )_{ik}^{xz} (\Phiml{\beta_1})_{kj}^{zy}a_{j,n+1}^y $$
$$= \sum\limits_{j=1}^{n} \sum\limits_{y \in \Z} (\Phiml{\beta_2}(\Phi_{\beta_1} )\Phiml{\beta_1})_{ij}^{xy}a_{j,n+1}^y $$

On the other hand, by definition, $\Phim{\beta_1\beta_2}(a_{i,n+1}^x) =  \sum\limits_{j=1}^{n} \sum\limits_{z \in \Z} (\Phiml{\beta_1\beta_2})_{ij}^{xy}a_{j,n+1}^y$.

Therefore, we have $(\Phiml{\beta_2}(\Phi_{\beta_1} )\Phiml{\beta_1})_{ij}^{xy} = (\Phiml{\beta_1\beta_2})_{ij}^{xy}.$
\end{proof}
\end{prop}

Let $I_n \in M_n(M_{\infty}(\A_n))$ be the identity matrix, i.e $(I_n)_{ij}^{xy} = \delta_{i,j}\delta_{x,y}.$ Then apparently, for a trivial braid $\beta \in \C_n, $ $\Phiml{\beta}, \Phimr{\beta}, \Phipl{\beta},\Phipr{\beta}$ are all equal to $I_n$. Therefore, we have the following corollary.

\begin{cor} \label{philr invertible}
For any braid $\beta \in \C_n,$ $\Phiml{\beta}, \Phimr{\beta}, \Phipl{\beta},\Phipr{\beta}$ are all invertible. Explicitly,

$$
(\Phiml{\beta})^{-1} = \Phiml{\beta^{-1}}(\Phi_{\beta} ), \quad
(\Phimr{\beta})^{-1} = \Phimr{\beta^{-1}}(\Phi_{\beta} ),
$$
$$
(\Phipl{\beta})^{-1} = \Phipl{\beta^{-1}}(\Phi_{\beta} ),\quad
(\Phipr{\beta})^{-1} = \Phipr{\beta^{-1}}(\Phi_{\beta} ).
$$
\begin{proof}
In Proposition \ref{Phi(beta 1 2)}, set $\beta_1 = \beta, \beta_2 = \beta^{-1}$.
\end{proof}
\end{cor}

\begin{prop} \label{key prop}
For any $\beta \in  \C_n,$ we have $\Phi_{\beta}(A) = \Phiml{\beta}A\Phimr{\beta} = \Phipl{\beta}A\Phipr{\beta}.$
\begin{proof}
By Proposition \ref{Phi(beta 1 2)}, it suffices to show the above equation holds for any $\alpha_k \in \C_n.$ This can be verified directly, though maybe tediously.

Here we provide another way to prove it.

$\Phi_{\beta}(A_{ij}^{xy}) = \Phi_{\beta}(a_{ij}^{x+y}) = \Phi_{\beta}(a_{i,n+1}^{x}*a_{n+1,j}^{y}) \\ \xlongequal{Proposition \, \ref{connectable}} \Phim{\beta}(a_{i,n+1}^{x})*\Phim{\beta}(a_{n+1,j}^{y}) \\ = (\sum\limits_{k=1}^{n} \sum\limits_{z \in \Z} (\Phiml{\beta})_{ik}^{xz}a_{k,n+1}^z) * (\sum\limits_{k'=1}^{n} \sum\limits_{z' \in \Z} a_{n+1,k'}^{z'}(\Phimr{\beta})_{k'j}^{z'y}) \\ = \sum\limits_{k,k'=1}^{n} \sum\limits_{z,z' \in \Z} (\Phiml{\beta})_{ik}^{xz}a_{kk'}^{z+z'}(\Phimr{\beta})_{k'j}^{z'y} \\ = \sum\limits_{k,k'=1}^{n} \sum\limits_{z,z' \in \Z} (\Phiml{\beta})_{ik}^{xz}A_{kk'}^{zz'}(\Phimr{\beta})_{k'j}^{z'y} = (\Phiml{\beta}A\Phimr{\beta})_{ij}^{xy}$

The other equation can be proved analogously.
\end{proof}
\end{prop}

\begin{cor} \label{A-PhiA}
For $\beta \in \C_n, 1 \leq p,q \leq n, f \in \Z,$ the elements of $A - \Lambda_{f;p,q}\Phi_{\beta}(A)\Lambda_{f;p,q}^{-1}$ are in $\mathcal{I}_{\beta;f;p,q}$. More generally, if $b = a_{i_1,i_2}^{x_1}a_{i_2,i_3}^{x_2} \cdots a_{i_k,i_{k+1}}^{x_k}, c_i = \lambda^{\delta_{i,p}}\mu^{-f\delta_{i,q}}$, then $b - c_{i_1}\Phi_{\beta}(b)c_{i_{k+1}}^{-1}$ is in $\mathcal{I}_{\beta;f;p,q}$.
\begin{proof}
Set $\Lambda = \Lambda_{f;p,q}.$ Then

$A - \Lambda\Phi_{\beta}(A)\Lambda^{-1} = A - \Lambda\Phiml{\beta}A\Phimr{\beta}\Lambda^{-1} = A - \Lambda\Phiml{\beta}A + \Lambda\Phiml{\beta}(A - A\Phimr{\beta}\Lambda^{-1}).$

The elements of the right hand side of the above equation are in $\mathcal{I}_{\beta;f;p,q}$, which implies the first part of the corollary. The more general statement in the corollary is then a direct consequence.
\end{proof}
\end{cor}

\subsection{Invariance proof}
\label{subsec:invariace proof}

In this subsection, we prove Theorem \ref{HC0 braid}. Apparently, the three parts in the theorem correspond to the three types of Markov moves introduced in Theorem \ref{markov move}. In the following three subsections, we prove each part of the theorem, respectively.

\subsubsection{Invariance under Markov Move I}
\label{subsubsec:markov I}

Let $\tilde{\beta} = \alpha^{-1}\beta\alpha, \, \alpha, \beta \in \C_n, f \in \Z$, and assume $\alpha(m) = 1$. Set $\Lambda_i = \Lambda_{f;i,i}$. We define an isomorphism $\varphi:$ $HC_0(\tilde{\beta};f; m) \longrightarrow HC_0(\beta;f;1)$ by specifying the image of the generators.
$$
\varphi(A):= \Phi_{\alpha}(A), \, i.e. \, \varphi(a_{ij}^x):= \Phi_{\alpha}(a_{ij}^x)
$$
We need to show $\varphi(\mathcal{I}_{\tilde{\beta};f;m,m}) \subset \mathcal{I}_{\beta;f;1,1}$.

First of all, by using Proposition \ref{Phi(beta 1 2)}, we have

$\Phi_{\alpha}(\Phiml{\alpha^{-1}\beta\alpha}) =\Phi_{\alpha}(\Phiml{\beta\alpha}(\Phi_{\alpha^{-1}}  )\Phiml{\alpha^{-1}}) = \Phiml{\beta\alpha}\Phiml{\alpha^{-1}}(\Phi_{\alpha}  ) = \Phiml{\alpha}(\Phi_{\beta}  )\Phiml{\beta}\Phiml{\alpha^{-1}}(\Phi_{\alpha}  ),$

Therefore, we have

\begin{align*}
 \varphi(A- \Lambda_m\Phiml{\tilde{\beta}}A) = & \varphi(A)- \Lambda_m\varphi(\Phiml{\tilde{\beta}})\varphi(A) \\
  =                                            & \Phi_{\alpha}(A) -\Lambda_m\Phi_{\alpha}(\Phiml{\alpha^{-1}\beta\alpha})\Phi_{\alpha}(A) \\
  =       & \Phi_{\alpha}(A) -\Lambda_m\Phiml{\alpha}(\Phi_{\beta})\Phiml{\beta}\Phiml{\alpha^{-1}}(\Phi_{\alpha})\Phiml{\alpha}A\Phimr{\alpha} \\
  =       & \Phiml{\alpha}A\Phimr{\alpha}- \Lambda_m\Phiml{\alpha}(\Phi_{\beta})\Phiml{\beta}A\Phimr{\alpha} \\
  =   & (\Phiml{\alpha}-\Lambda_m\Phiml{\alpha}(\Phi_{\beta})\Lambda_1^{-1})A\Phimr{\alpha} +  \\
      &  \Lambda_m\Phiml{\alpha}(\Phi_{\beta})\Lambda_1^{-1}(A - \Lambda_{1}\Phiml{\beta}A)\Phimr{\alpha}\\
\end{align*}

Since $(\Phiml{\alpha})_{ij}^{xy}$ is a non-commutative polynomial in which each monomial is of the form $a_{\alpha(i),j_1}^{x_1}a_{j_1,j_2}^{x_2} \cdots a_{j_{k-1},j}^{x_k}$, and note that $\delta_{i,m} = \delta_{\alpha(i),1}$, then $(\Phiml{\alpha}-\Lambda_m\Phiml{\alpha}(\Phi_{\beta})\Lambda_1^{-1})_{ij}^{xy}$ is a sum of polynomials of the form $a_{\alpha(i),j_1}^{x_1}a_{j_1,j_2}^{x_2} \cdots a_{j_{k-1},j}^{x_k} - (\lambda\mu^{-f})^{\delta_{\alpha(i),1}}\Phi_{\beta}(a_{\alpha(i),j_1}^{x_1}a_{j_1,j_2}^{x_2} \cdots a_{j_{k-1},j}^{x_k})(\lambda\mu^{-f})^{-\delta_{j,1}}$, which, by Corollary \ref{A-PhiA}, is in $\mathcal{I}_{\beta;f;1,1}$.

Since elements of $A-\Lambda_1\Phiml{\beta}A$ are also in $\mathcal{I}_{\beta;f;1,1}$, this implies $\varphi(A- \Lambda_m\Phiml{\tilde{\beta}}A) \subset \mathcal{I}_{\beta;f;1,1}$.

%
%

The proofs of the other three relations $ \varphi(A- A\Phimr{\tilde{\beta}}\Lambda_m^{-1}), \varphi(A- \Lambda_m\Phipl{\tilde{\beta}}A),\,  \varphi(A- A\Phipr{\tilde{\beta}}\Lambda_m^{-1})$ are completely analogous.

This shows $\varphi(\mathcal{I}_{\tilde{\beta}};f;m,m) \subset \varphi(\mathcal{I}_{\beta};f;1,1)$ and thus induces a well-defined map $HC_0(\tilde{\beta};f) \longrightarrow HC_0(\beta;f)$. In a similar way, we can define the inverse map $HC_0(\beta;f) \longrightarrow HC_0(\tilde{\beta};f)$ by sending $A$ to $\Phi_{\alpha^{-1}}(A)$ and show that it is well defined. Thus $\varphi$ is an isomorphism.

\subsubsection{Invariance under Markov Move II}
\label{subsubsec:markov II}

For any $\beta \in \C_n, f \in \Z$ let $\tilde{\beta} = \Em(\beta)\alpha_{n}$. We show $HC_0(\tilde{\beta};f) \simeq HC_0(\beta;f+1).$

\begin{rem}
The proof of $HC_0(\Em(\beta)\alpha_{n}^{-1};f+1) \simeq HC_0(\beta;f)$ is completely analogous. To save space, we omit its proof here.
\end{rem}

Define $\varphi: HC_0(\tilde{\beta};f;n+1) \longrightarrow HC_0(\beta;f+1;n),$

\begin{equation}
\varphi(a_{ij}^x) =
\begin{cases}
a_{nn}^x & i = n+1, j = n+1\\
\mu a_{nj}^x & i =n+1, j \leq n \\
\mu^{-1}a_{in}^x & i\leq n, j = n+1 \\
a_{ij}^x & i\leq n, j \leq n\\
\end{cases}
\end{equation}

The verification that $\varphi$ maps $\mathcal{I}_{\tilde{\beta};f;n+1,n+1}$ to $\mathcal{I}_{\beta;f+1;n,n}$ consists of direct but long calculations. we will only show $\varphi(a_{i,j}^{x+y} - (\lambda\mu^{-f})^{\delta_{i,n+1}} \Phim{\tilde{\beta}}(a_{i,n+2}^x)*a_{n+2,j}^y) \in \mathcal{I}_{\beta;f+1;n,n}$. The other relations can be proven similarly.

Set $c = \lambda\mu^{-f}, \mathcal{I} = \mathcal{I}_{\beta;f+1;n,n}$.

$\\$
Case $1: i=n+1, j\leq n.$

$\varphi(a_{n+1,j}^{x+y} - c \Phim{\tilde{\beta}}(a_{n+1,n+2}^x)*a_{n+2,j}^y)
= \mu a_{n,j}^{x+y} - c \varphi(\Phim{\Em(\beta)}(a_{n,n+2}^x)*a_{n+2,j}^y) \\
= \mu a_{n,j}^{x+y} - c \Phim{\beta}(a_{n,n+1}^x)*a_{n+1,j}^y
= \mu(a_{n,j}^{x+y} - \lambda\mu^{-f-1} \Phim{\beta}(a_{n,n+1}^x)*a_{n+1,j}^y) \in \mathcal{I}$

Note that here we used the fact that $\Phim{\Em(\beta)}(a_{n,n+2}^x)*a_{n+2,j}^y = \Phim{\beta}(a_{n,n+1}^x)*a_{n+1,j}^y \in \A_{n+1}$.
$\\$
Case $2: i= n, j \leq n.$

$a_{n,j}^{x+y} - \Phim{\tilde{\beta}}(a_{n,n+2}^x)*a_{n+2,j}^y \\
= a_{n,j}^{x+y} - \Phim{\Em(\beta)}(-a_{n+1,n+2}^x + \frac{1}{\Var\mu}a_{n+1,n}^0 a_{n,n+2}^x)*a_{n+2,j}^y \\
= a_{n,j}^{x+y} - (-a_{n+1,n+2}^x + \frac{1}{\Var\mu}\Phim{\Em(\beta)}(a_{n+1,n}^0)\Phim{\Em(\beta)}(a_{n,n+2}^x))*a_{n+2,j}^y \\
=  a_{n,j}^{x+y} + a_{n+1,j}^{x+y} - \frac{1}{\Var\mu}\Phi_{\Em(\beta)}(a_{n+1,n}^0)\Phim{\Em(\beta)}(a_{n,n+2}^x)*a_{n+2,j}^y \\
=  a_{n,j}^{x+y} + a_{n+1,j}^{x+y} - \frac{1}{\Var\mu}\Phim{\beta}(a_{n+1,n}^0)\Phim{\beta}(a_{n,n+1}^x)*a_{n+1,j}^y$

Since $\varphi(\Phim{\beta}(a_{n+1,n}^0)) = \mu a_{n,n+1}^0 * \Phim{\beta}(a_{n+1,n}^0) = c a_{nn}^{0} \; (\; \textrm{mod}\; \mathcal{I}),$

$\varphi(a_{n,j}^{x+y} - \Phim{\tilde{\beta}}(a_{n,n+2}^x)*a_{n+2,j}^y) \\
= (1+u)a_{nj}^{x+y} - \frac{1}{\Var\mu} c (1+u)\Var \Phim{\beta}(a_{n,n+1}^x)*a_{n+1,j}^y \;(\; \textrm{mod}\; \mathcal{I}) \\
= (1+u)(a_{nj}^{x+y} - \lambda\mu^{-f-1}\Phim{\beta}(a_{n,n+1}^x)*a_{n+1,j}^y )(\; \textrm{mod}\; \mathcal{I}) = 0\;(\; \textrm{mod}\; \mathcal{I}).$
$\\$
Case $3: i\leq n-1, j \leq n$.

$\varphi(a_{i,j}^{x+y} - \Phim{\tilde{\beta}}(a_{i,n+2}^x)*a_{n+2,j}^y) = \varphi(a_{i,j}^{x+y} - \Phim{\Em(\beta)}(a_{i,n+2}^x)*a_{n+2,j}^y)\\
 = \varphi(a_{i,j}^{x+y} - \Phim{\beta}(a_{i,n+1}^x)*a_{n+1,j}^y) = a_{i,j}^{x+y} - \Phim{\beta}(a_{i,n+1}^x)*a_{n+1,j}^y \in \mathcal{I}.$
$\\$
Case $4: j = n+1.$ The proof is the same as the above threes cases except an overall scalar $\mu^{-1}$ is multiplied to each expression.

This finishes the verification. One can also define a map $HC_0(\beta;f+1;n) \longrightarrow HC_0(\tilde{\beta};f;n+1) $ sending $a_{ij}^x$ to $a_{ij}^x$, and show that it is well defined. Clearly this is the inverse of $\varphi$.

\subsubsection{Invariance under Markov move III} \label{invariance III}

Recall that $D_n$ is the unit disk with $n+1$ punctures $p, p_1, \cdots, p_n$ centered at the origin of the complex plane. To be more precise, let $p$ be the origin and the coordinate of $p_i$ be $\frac{i}{n+1}.$ We define a map $r: D_n \longrightarrow D_n$ by $r(z) = \frac{\bar{z}}{|z|} - \bar{z}.$ Namely, $r$ is a reflection about the $x$-axis followed by another reflection about the circle centered at the origin with radius $\frac{1}{2}$. Note that $r^2 = Id$. Also $r \times Id$ defines a map on $X = D_n \times [0,1]$, which will still be denoted by $r$.

Recall that $\C_n$ is the braid group on the punctured disk $D_n$ inside $X$. Therefore, $r$ induces a group isomorphism from $\C_n$ to itself. Explicitly, the isomorphism, also denote by $r$, is given by:

\begin{equation}
r(\alpha_i) =
\begin{cases}
(\alpha_{n-1} \cdots \alpha_1 \alpha_0 \alpha_1 \cdots \alpha_{n-1})^{-1} & i = 0 \\
\alpha_{n-i}                                                  & 1 \leq i \leq n-1 \\
\end{cases}
\end{equation}

\begin{lem}
The map $r$ defined above from $\C_n$ to $\C_n$ is a group isomorphism and $r^2 = Id$.
\begin{proof}
This can be verified purely algebraically.
\end{proof}
\end{lem}

Also recall that $q_1,  \cdots, q_n$ are $n$ points with the coordinate $\frac{i}{n+1} -\epsilon$ for some tiny $\epsilon>0$. And $Q_n = \{ q_i, 1 \leq i \leq n\}$, $\Q_n = \{\gamma:[0,1] \longrightarrow D_n | \, \gamma  \\  \textrm{is continuous}, \gamma(0),\gamma(1) \in Q_n\} / \thicksim$. Let $q'_{n+1-i} = r(q_i),$ which has the coordinate $\frac{n+1-i}{n+1}+ \epsilon$, and let $Q_n' = \{ q_i', 1 \leq i \leq n\}$. It should be clear that in the definition of $\tilde{\A}_{n}$, if we replace $q_i$ by $q_i'$, insist that the curves start and end at $q_i'$, and change the skein relations accordingly, then we get the same algebra.

For a curve $\gamma \in \Q_n$ from $q_i$ to $q_j$,  $r(\gamma) $ is a curve from $q_{n+1-i}'$ to $q_{n+1-j}'$. And it's not hard to check that this map also preserves the \lq\lq skein" relations in Figure \ref{skein12} that defines $\tilde{A}_n.$ Thus $r$ induces an algebra isomorphism from $\tilde{A}_n$ to $\tilde{A}_n$.

\begin{figure}[h!]
\def\svgwidth{9cm}
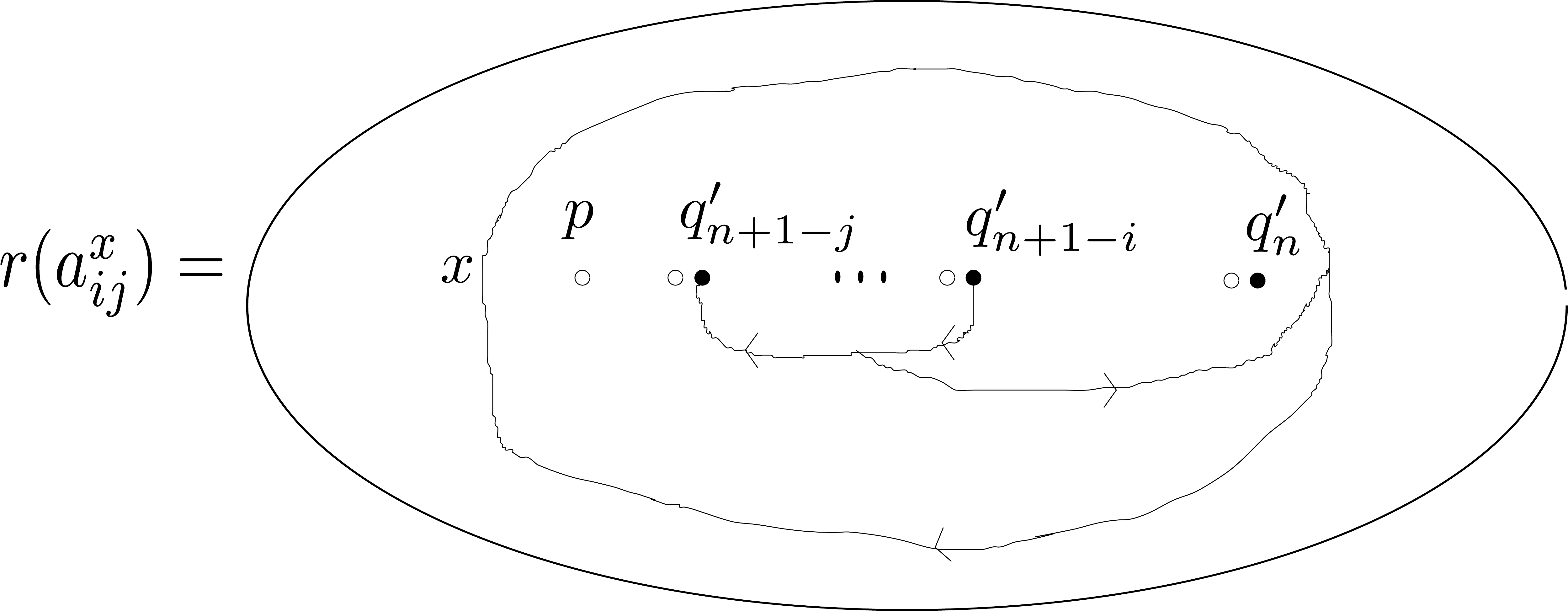
\caption{$r(a_{ij}^x)$}
\label{r_gammaij}
\end{figure}

Explicitly, the map $r: \tilde{A}_n \longrightarrow \tilde{A}_n$ is given by Figure \ref{r_gammaij}.

\begin{rem}
$r$ also extends to a bijection from $\Q_n^{+}$ to $\Q_n^{-}$ by furthermore requiring that $p_0$ is mapped to $p_{n+1}$. And $r$ maps the \lq\lq skein" relations that define $\Ap$ to the corresponding \lq\lq skein" relations that define $\Am$. Consequently, we get an isomorphism $r: \Ap \longrightarrow \Am.$ Note that the inverse map is also induced by $r$ that maps $\Q_n^{-}$ to $\Q_n^{+}$. For this reason, we will denote the inverse map also by $r$. In summary, $r$ is an isomorphism between $\Ap$ and $\Am$, which restricts to an isomorphism on $\An$ and which has square $Id$.
\end{rem}

\begin{lem} \label{star and r}
If $P,Q \in \Apm$ are connectable, then $r(P)$, $r(Q)$ are connectable, and $r(P*Q) = r(P)*r(Q)$.
\begin{proof}
%
Clear from the geometrical interpretation of $a_{ij}^x$ and the map $r$.
\end{proof}
\end{lem}

\begin{lem} \label{Phi and r}
If $\beta$ is a braid in $\C_n,$ then we have $r \circ \Phi_{\beta} = \Phi_{r(\beta)} \circ r$. More generally, we have $r \circ \Phim{\beta} = \Phip{r(\beta)} \circ r$
\begin{proof}

It's possible, though tedious, to prove it algebraically. For example, it suffices to prove the case for $\beta = \alpha_{k}^{\pm 1} $ acting on $a_{ij}^x$. Here we give another geometric proof which makes the statement in the lemma almost trivial. Recall that the isomorphism $r : \C_n \longrightarrow \C_n$ is induced by the homeomorphism $r \times Id : D_n \times I \longrightarrow D_n \times I$. By Remark \ref{rem picture action}, $\beta(\gamma_{ij}^x)$ can be obtained as the curve by sliding $\gamma_{ij}^x$ in $D_n \times \{0\}$ along the parallel copy braid diagram $\beta'$ up to $D_n \times \{1\}$. The map $r \times Id$ maps $\gamma_{ij}^x$ to $r(\gamma_{ij}^x)$, $\Phi_{\beta}(\gamma_{ij}^x)$ to $r \circ \Phi_{\beta}(\gamma_{ij}^x)$, and $\beta$ to $r(\beta)$. Thus $r \circ \Phi_{\beta}(\gamma_{ij}^x)$ is obtained by sliding $r(\gamma_{ij}^x)$ along the parallel copy braid diagram  $r(\beta)'$, and therefore $\Phi_{r(\beta)} \circ r (\gamma_{ij}^x) = r \circ \Phi_{\beta}(\gamma_{ij}^x).$

The more general equation can be proved analogously by using Remark \ref{rem picture action} and Lemma \ref{lem epsilonpm}.
\end{proof}
\end{lem}

\begin{thm} \label{r isomorphism}
For $\beta \in \C_n, f \in \Z$ the map $r: \A_n \longrightarrow \A_n$ induces an isomorphism from $HC_0(\beta;f;1) $ to $HC_0(r(\beta);f;n).$
\begin{proof}
It suffices to show $r$ maps $I_{\beta;f;1,1}$ to $I_{r(\beta);f;n,n}.$ Set $c = \lambda\mu^{-f}$.

$r((\Lambda_{f;1,1}\Phiml{\beta}A)_{ij}^{xy}) =  r(c^{\delta_{i,1}}\Phim{\beta}(a_{i,n+1}^x) * a_{n+1,j}^y) =  c^{\delta_{i,1}} (r \circ \Phim{\beta}(a_{i,n+1}^x)) * r(a_{n+1,j}^y) = c^{\delta_{n+1-i,n}} (\Phip{r(\beta)} \circ r (a_{i,n+1}^x)) * r(a_{n+1,j}^y).$

The first identity in the above equation is by the argument in Part $(2)$ of Remark \ref{rem relation}, the second identity is by Lemma \ref{star and r}, and the third by Lemma \ref{Phi and r}.

Assume $r(a_{i,n+1}^x) = \sum P_{k}^z a_{k0}^z, \, r(a_{n+1,j}^y) = \sum a_{0k'}^{z'}Q_{k'}^{z'},$ where $P_{k}^z, Q_{k'}^{z'}$ are elements in $\A_n$. Then

$r((A-\Lambda_{f;1,1}\Phiml{\beta}A)_{ij}^{xy})\\
= \sum P_{k}^z a_{kk'}^{zz'}Q_{k'}^{z'} -  c^{\delta_{n+1-i,n}}\Phi_{r(\beta)}(P_{k}^z)\Phip{r(\beta)}(a_{k0}^z) * a_{0k'}^{z'}Q_{k'}^{z'} \\
= \sum(P_{k}^z - c^{\delta_{n+1-i,n}}\Phi_{r(\beta)}(P_{k}^z)c^{-\delta_{k,n}})a_{kk'}^{zz'}Q_{k'}^{z'} + c^{\delta_{n+1-i,n}}\Phi_{r(\beta)}(P_{k}^z)c^{-\delta_{k,n}}(a_{kk'}^{zz'} - c^{\delta_{k,n}}\Phip{r(\beta)}(a_{k0}^z) * a_{0k'}^{z'})Q_{k'}^{z'}$

Note that $P_k^z$ is a sum of monomials of the form $a_{n+1-i,i_1}^{x_1}a_{i_1,i_2}^{x_2}\cdots a_{i_{m-1},k}^{x_m}$, then $P_{k}^z - c^{\delta_{n+1-i,n}}\Phi_{r(\beta)}(P_{k}^z)c^{-\delta_{k,n}}$ is in $I_{r(\beta);f;n,n}$ by Corollary \ref{A-PhiA}.

Then it follows that $r((A-\Lambda_{f;1,1}\Phiml{\beta}A)_{ij}^{xy})$ is in $I_{r(\beta);f;n,n}$.

The other relations are proved in basically the same way. And thus we showed $r$ is well-defined. The fact that $r$ is an isomorphism is trivial to check.

\end{proof}
\end{thm}

Now we prove $HC_0(\beta;f)$ is invariant under Markov move III. A key observation is the following commuting diagram.
\begin{equation}
\xymatrix{
\C_n \ar[r]^-{\Ep} \ar[d]^-{r} & \C_{n+1} \ar[d]^-{r} \\
\C_n \ar[r]^-{\Em}            & \C_{n+1}\\
}
\end{equation}

\begin{lem} \label{lem epsilonpm}
The above diagram commutes, namely $r \circ \Ep = \Em \circ r : \C_n \longrightarrow \C_{n+1}$.
\begin{proof}
We only need to check on the generators.

$r \Ep(\alpha_0) = r(\alpha_1\alpha_0\alpha_1) = \alpha_{n} (\alpha_{n} \cdots \alpha_1 \alpha_0 \alpha_1 \cdots \alpha_{n})^{-1}\alpha_{n} \\ = (\alpha_{n-1} \cdots \alpha_1 \alpha_0 \alpha_1 \cdots \alpha_{n-1})^{-1} = \Em r(\alpha_0).$

For $i \geq 1,$ $r \Ep(\alpha_i) = r(\alpha_{i+1}) = \alpha_{n-i} = \Em r(\alpha_i)$
\end{proof}
\end{lem}

Let $\beta \in \C_n, f \in \Z$, then $r(\Ep(\beta)\alpha_1^{\pm 1}) = r(\Ep(\beta))r(\alpha_1^{\pm 1}) = \Em(r(\beta))\alpha_{n}^{ \pm 1}.$ Therefore,

$HC_0(\Ep(\beta)\alpha_1^{\pm 1};f) \simeq HC_0(r(\Ep(\beta)\alpha_1^{\pm 1});f) = HC_0(\Em(r(\beta))\alpha_{n}^{ \pm 1};f) \\ \simeq HC_0(r(\beta); f \pm 1) \simeq HC_0(\beta;f \pm 1).$

The first and last isomorphism above are due to Theorem \ref{r isomorphism} and the second isomorphism is the invariance isomorphism under Markov move II.

Now we finished showing $HC_0(\beta;f)$ is invariant under Markov move III.

\section{Properties of the invariant}

\subsection{Symmetries of the invariant} \label{subsec:symmetry}

In Section \ref{invariance III}, we already proved that for a braid $\beta \in \C_n$, we have $HC_0(\beta;f) \simeq HC_0(r(\beta);f)$. Here we show the relation between $HC_0(\beta;f)$ and $HC_0(\beta^{-1};f)$.

\begin{prop}
Let $\beta \in \C_n, f \in \Z$, then $HC_0(\beta^{-1};f) $ is isomorphic to $HC_0(\beta;-f)$ with $\lambda$ replaced by $\lambda^{-1}$.
\begin{proof}
Let $HC_0'(\beta;-f)$ be the algebra obtained from $HC_0(\beta;-f)$ by replacing $\lambda$ by $\lambda^{-1}$. We define the isomorphism $HC_0(\beta^{-1};f) \longrightarrow HC_0'(\beta;-f)$ to be the one induced by $\Phi_{\beta}$. We need to check $\Phi_{\beta}$ maps $\mathcal{I}_{\beta^{-1};f;1,1}$ to $\mathcal{I}_{\beta;-f;1,1}$ with $\lambda$ replaced by $\lambda^{-1}$. Set $\Lambda = \Lambda_{f;1,1}$, and note that $\Lambda^{-1}$ is exactly the matrix $\Lambda_{-f;1,1}$ with $\lambda$ replaced by $\lambda^{-1}$.

$\Phi_{\beta}(\Lambda\Phipl{\beta^{-1}}A - A) = \Lambda\Phipl{\beta^{-1}}(\Phi_{\beta}) \Phi_{\beta}(A) - \Phi_{\beta}(A) = \Lambda\Phipl{\beta^{-1}}(\Phi_{\beta})\Phipl{\beta}A\Phipr{\beta} - \Phipl{\beta}A\Phipr{\beta} = \Lambda A\Phipr{\beta} - \Phipl{\beta}A\Phipr{\beta} = \Lambda(A - \Lambda^{-1}\Phipl{\beta}A)\Phipr{\beta} $

The second equality is by Proposition \ref{key prop} and the third one is by Corollary \ref{philr invertible}.

The other three relations can be proved analogously. Therefore, $\Phi_{\beta}$ induces a well-defined algebra map from $HC_0(\beta^{-1};f) $ to $HC_0'(\beta;-f)$. It's easy to check it's also an isomorphism.
\end{proof}
\end{prop}

\subsection{Torus Knots} \label{subsec:torus knot}

In this subsection we study some properties of the torus knots in $\Lens$.

Let $C$ be the equator of $S^2$, then $S^1 \times C$ is a torus which bounds two solid tori in $\Lens$, with $z_0 \times C$ being the meridian and $S^1 \times z_1$ the longitude. In \cite{Ding}, a knot in $\Lens$ is called a torus knot if it can be isotoped to a knot in $S^1 \times C$. Fix a meridian $M$ and a longitude $L$ in $S^1 \times C$, and let $p,q$ be two relatively prime integers. A $(p,q)$-knot in $\Lens$ is a knot which can be isotoped to $p M + q L$ in $S^1 \times C$. In general, for a knot $K$ and a framing $l$, $HC_0(K;l)$ may not be finitely generated as an $R$-algebra. However, we show below that for torus knots, the invariant indeed is always finitely generated.

\begin{thm} \label{thm torus knot}
Let $K$ be a $(p,q)$-knot in $\Lens$ with framing $l$ where $p,q$ are relatively prime integers, then $HC_0(K;l)$ is finitely generated as an $R$-algebra. Moreover, the minimum number of algebra generators is no more than $q-1$.

\begin{proof}
By Remark \ref{S1S2 model}, a $(p,q)$-knot is represented by the braid $\beta(p,q) = (\alpha_0 \cdots \alpha_{p-1})^q$. See Figure \ref{torus knot} for a picture of $(3,2)$-knot. For simplicity, we still use $\beta$ to denote $\beta(p,q)$. Also for reasons that will become clear below, we use the notation $b_{ij}^x = a_{i+1,j+1}^x$. Assume $HC_0(K;l) = HC_0(\beta;f) = \A_{p} / \mathcal{I}_{\beta;f;1,1}$, and set $c = \lambda\mu^{-f}$. It's easy to check that the following equation holds:

\begin{figure}[h!]
\centering
\def\svgwidth{3cm}
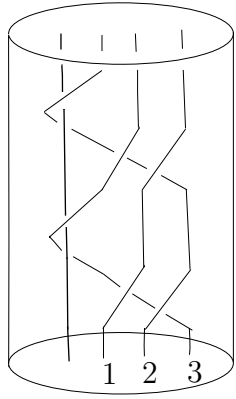
\caption{$(3,2)$-knot}
 \label{torus knot}
\end{figure}

\begin{equation}
\Phip{\beta(p,1)}(a_{i0}^x) =
\begin{cases}
a_{i+1,0}^x &  1 \leq i \leq p-1 \\
\mu a_{1,0}^{x-1}  &  i = p \\
\end{cases}
\end{equation}

Then we have $\Phip{\beta(p,q)}(a_{i0}^x) = \mu^{\lfloor \frac{i-1+q}{p}\rfloor}a_{(i-1+q) \, (mod \, p) + 1, 0}^{x - \lfloor \frac{i-1+q}{p}\rfloor}$. Using $b_{ij}^x$ to replace $a_{i+1,j+1}^x$, we get a simpler expression $\Phip{\beta(p,q)}(b_{i,-1}^x) = \mu^{\lfloor \frac{i+q}{p}\rfloor}b_{(i+q) \, (mod \, p), -1}^{x - \lfloor \frac{i+q}{p}\rfloor}.$

Thus by Part $2$ of Remark \ref{rem relation}, the third relation that defines $\mathcal{I}_{\beta;f;1,1}$ is

\begin{equation} \label{b_ij i}
b_{ij}^{x} - \mu^{\lfloor \frac{i+q}{p}\rfloor} c^{\delta_{i,0}} b_{(i+q) \, (mod \, p), j}^{x - \lfloor \frac{i+q}{p}\rfloor} , \, 0 \leq i,j \leq p-1, x \in \Z.
\end{equation}

Similarly, the fourth relation that defines $\mathcal{I}_{\beta;f;1,1}$ is

\begin{equation} \label{b_ij j}
b_{ij}^{x} - \mu^{-\lfloor \frac{j+q}{p}\rfloor} c^{-\delta_{j,0}}b_{i, (j+q) \, (mod \, p)}^{x + \lfloor \frac{j+q}{p}\rfloor} , \, 0 \leq i,j \leq p-1, x \in \Z.
\end{equation}

Define $g(i,k):= \sum\limits_{r=0}^{k-1} \lfloor \frac{(i+rq) \, (mod \, p) + q}{p} \rfloor,\, h(i,k):= \sum\limits_{r=0}^{k-1} \delta_{(i+rq) \, (mod \, p),0}, \,   0 \leq i \leq p-1, k \geq 1$, and define $f(i,0):= 0, h(i,0):= 0$.

It's elementary to check that $g(i,k) = \lfloor\frac{k}{p}\rfloor q + g(i, k \, mod \, p)$ and $h(i,k) = \lfloor\frac{k}{p}\rfloor + h(i, k \, mod \, p)$, and in $HC_0(\beta;f;1,1)$, we have the equalities $b_{ij}^x = \mu^{g(i,k)}c^{h(i,k)}b_{(i + kq) \, (mod \, p),j}^{x - g(i,k)} = \mu^{-g(j,k)}c^{-h(j,k)}b_{i,(j + kq) \, (mod \, p)}^{x + g(j,k)}, \forall k \geq 0.$ Especially, we have $b_{ij}^x = \mu^{g(i,p)}c^{h(i,p)}b_{ij}^{x - g(i,p)} = \mu^{q}c b_{ij}^{x-q},$ so $b_{ij}^x$ is periodic, up to a scalar, in $x$ with period equal to $q$.

Let $k_1,k_2$ be any numbers that satisfy $k_1 q \, (mod \, p) = i, \, k_2 q \, (mod \, p) = j$, then $b_{ij}^x = \mu^{g(0,k_2) - g(0,k_1)}c^{h(0,k_2) - h(0,k_1)}b_{00}^{x + g(0,k_1) - g(0,k_2)},$ and $b_{00}^{x+q} = \mu^{q}c b_{00}^{x}.$ Thus all the $b_{ij}^x \;'$s are completely determined by $b_{00}^0 = (1+\mu)\Var, b_{00}^1, \cdots, b_{00}^{q-1}$ and the condition that $b_{00}^{x+q} = \mu^{q}c b_{00}^{x}.$ So $HC_0(\beta;f)$ is finitely generated and $\{b_{00}^x,  1 \leq x \leq q-1\}$ is a set of generators.

\end{proof}
\end{thm}

At the end of this subsection, let's compute some examples of torus knots.
\begin{example}

$1).$ \textbf{\bm{$(p,1)$}-knot}. The $(p,1)$-knot is represented by the braid $\alpha_0 \cdots \alpha_{p-1}$. Clearly, by Markov II in Theorem \ref{markov move}, this braid is equivalent to $\alpha_0$ representing the $(1,1)$-knot. Set $\beta = \alpha_0 \in \C_1,\Lambda = \Lambda_{0;1,1}, f = 0$. By Theorem \ref{thm torus knot}, $a_{11}^{x+1} = \lambda\mu a_{11}^x$. Since $a_{11}^0 = (1+\mu)\Var,$ we have $a_{11}^x = (1+\mu)\Var (\lambda\mu)^x$.

By definition, $\Phim{\beta}(a_{12}^x) = -\mu a_{12}^{x-1} + \frac{1}{\Var}a_{11}^{x}a_{12}^{-1},$ thus $(\Lambda\Phiml{\beta}A-A)_{11}^{xy} = -\lambda\mu a_{11}^{x+y-1} + \frac{\lambda}{\Var}a_{11}^{x}a_{12}^{y-1} - a_{11}^{x+y}$. The second relation can be calculated analogously. By using the fact that $a_{11}^x = (1+\mu)\Var (\lambda\mu)^x$, it's easy to see that $HC_0(\alpha_0) \simeq R / \langle \mu^2 -1\rangle$.

$2).$ \textbf{\bm{$(p,2)$}-knot}. By Proposition $2.2$ in \cite{Ding}, all the $(p,2)$-knots are equivalent to each other with $p$ odd.
This can also been seen directly by Markov moves.
Thus, we only need to compute the $(1,2)$-knot, which is represented by $\beta = \alpha_0^2$. It was shown in the second example in Section \ref{subsec:example} that $HC_0(\alpha_0^2) \simeq R[X]/\langle (1-\mu)X, X^2 - \Var^2\lambda(1+\mu)^2\rangle$.

$3).$ \textbf{\bm{$(p,3)$}-knot}. Again by Proposition $2.2$ in \cite{Ding}, there are two classes of knots of this type. A representative of each class could be chosen as $(1,3)$-knot and $(2,3)$-knot. Here we only compute $HC_0(\alpha_0^3)$. Since the calculations are not difficult but tedious, we just present the result obtained by computer packages. $HC_0(\alpha_0^3) \simeq R\langle X,Y\rangle / \langle Y^2 - \Var\lambda\mu^4(1+\mu^2)X, X^2 - \Var(1+\mu^{-2})Y, (1+\mu^2)(XY-YX), -\mu^2 XY + YX + \Var^2\lambda\mu^4(\mu^2-1) \rangle$.

\end{example}

\subsection{Local knots} \label{subsec:local knots}

Throughout this subsection, we will set $\Var = -1$. A knot is called {\it{local}} if it is contained in a $3$-ball. It's easy to see that a knot in $\Lens$ is local if and only if it can be represented as the closure of a braid which doesn't contain $\alpha_0$ or $\alpha_0^{-1}$, i.e. a braid in $\B_n = \langle \alpha_1, \cdots, \alpha_{n-1} \rangle \subset \C_n$. Note that the braids in $\B_n$ are closed under the Markov moves given in Theorem \ref{markov move}, and moreover, Markov move III in this case is a consequence of Markov moves I, II. Since Markov moves I, II are just the classical Markov moves for braids in $\B_n$ representing knots in $S^3$, we thus have a one-to-one correspondence between knots in $S^3$ and local knots in $\Lens$.

Let $\beta \in \B_n \subset \C_n$, and let $hc_0(\beta)$ denote the $0$-th framed knot contact homology in \cite{ng2008framed}. Then we have the following decompositions for the $HC_0$ invariant of local knots, which relates our invariant to $hc_0$.

\begin{prop} \label{HC0 decomposition}
Let $\beta \in \B_n$ be a braid such that its closure is a local knot in $\Lens$, then $HC_0(\beta) \simeq hc_0(\beta) + \sum\limits_{0 \neq x \in \Z} H_x$, where all the $H_x\;'$s are isomorphic to each other as subalgebras and there is a surjective algebra morphism from $H_x$ to $hc_0(\beta)$. Moreover, for any $0 < m \in \Z$, $\sum\limits_{|x|\leq m} H_x$ is a proper subalgebra of $HC_0(\beta)$.
\begin{proof}

Set $\Lambda = \Lambda_{0;1,1}$.

It's easy to see from Equation \ref{equ:Phi_k} that $\Phim{\beta}(a_{i, n+1}^x) = \Phim{\beta}(a_{i,n+1}^0) * a_{n+1,n+1}^x$ and that $\Phip{\beta}(a_{i,0}^x) = \Phim{\beta}(a_{i,n+1}^x) * a_{n+1,0}^0 = \Phim{\beta}(a_{i,n+1}^0) * a_{n+1,0}^x$. Therefore, $(\Phipl{\beta}A)_{ij}^{xy} = \Phip{\beta}(a_{i,0}^{x})*a_{0,j}^{y} = \Phim{\beta}(a_{i,n+1}^x)*a_{n+1,j}^y = (\Phiml{\beta}A)_{ij}^{xy}$, i.e. $\Phipl{\beta}A = \Phiml{\beta}A$. Similarly, we have $A\Phipr{\beta} = A\Phimr{\beta}$. So to compute $HC_0(\beta)$, we only need to consider the relations $A - \Lambda\Phiml{\beta}A, A - A\Phimr{\beta}\Lambda^{-1}.$

$(\Lambda\Phiml{\beta}A)_{ij}^{xy} = \lambda^{\delta_{i,1}}\Phim{\beta}(a_{i,n+1}^x)*a_{n+1,j}^y = \lambda^{\delta_{i,1}}\Phim{\beta}(a_{i,n+1}^0)*a_{n+1,j}^{x+y} = \lambda^{\delta_{i,1}}\Phim{\beta}(a_{i,n+1}^0)*a_{n+1,j}^{0}*a_{j,j}^{x+y} = (\Lambda\Phiml{\beta}A)_{ij}^{00}*a_{j,j}^{x+y}$.

Therefore, we have $A_{ij}^{xy} - (\Lambda\Phiml{\beta}A)_{ij}^{xy} = (A_{ij}^{00} - (\Lambda\Phiml{\beta}A)_{ij}^{00})*a_{jj}^{x+y}.$ Similarly, we have $A_{ij}^{xy} - (A\Phimr{\beta}\Lambda^{-1})_{ij}^{xy} = a_{ii}^{x+y}*(A_{ij}^{00} - (A\Phimr{\beta}\Lambda^{-1})_{ij}^{00}).$

Let $E_0 = \Z\langle a_{ij}^0, 1 \leq i,j \leq n \rangle$ and for $ x \neq 0, \, E_x = \Z\langle a_{ij}^0,a_{ij}^x,  1 \leq i,j \leq n \rangle$. For all $x \in \Z$, let $H_x = E_{x} / \langle (A_{ij}^{00} - (\Lambda\Phiml{\beta}A)_{ij}^{00})*a_{jj}^x, a_{ii}^x*(A_{ij}^{00} - (A\Phimr{\beta}\Lambda^{-1})_{ij}^{00}), 1 \leq i,j \leq n \rangle$. Then clearly for $x \neq 0, y \neq 0,$ $H_x$ is isomorphic to $H_y$ by sending $a_{ij}^0$ to $a_{ij}^0$ and $a_{ij}^x$ to $a_{ij}^y$. And there is also a surjective morphism from $H_x$ to $H_0$ by sending $a_{ij}^0, a_{ij}^x $ both to $a_{ij}^0$. Similarly, for all $x \in \Z$, one can define a surjective morphism from $HC_0(\beta)$ to $H_x$ by sending $a_{ij}^x$ to $a_{ij}^x$ and $a_{ij}^z$ to $a_{ij}^0$ for all $z \neq x$. We denote this surjection by $\pi_x$. It's also easy to see from the definition that there is an algebra morphism $\iota_x: H_x \longrightarrow HC_0(\beta)$ such that $\iota_x(a_{ij}^0) = a_{ij}^0, \iota_x(a_{ij}^x) = a_{ij}^x,$ and thus $\pi_x \iota_x = Id$. Therefore, we conclude that for all $x \in \Z$, $H_x$ is a subalgebra of $HC_0(\beta)$,  and that $HC_0(\beta) = H_0 + \sum\limits_{0 \neq x } H_x$, and that for all $m > 0$, $\sum\limits_{|x| \leq m } H_x$ is a proper subalgebra since it does not contain $H_y$ for $|y| > m$.

$\\$
Next we show that $H_0 \simeq hc_0(\beta)$.

It could be checked that in Equation \ref{equ:Phi_k}, if we set $\Var = -1, x = 0$, then $\Phi_{\beta}$ acting on $E_0$ is exactly the same as the braid action given in \cite{ng2008framed} if we make the change of variables as follows:  $a_{ij}^0 = \mu a_{ij}$ if $i> j$ and $a_{ij}^0 = a_{ij}$ otherwise. Note that here $a_{ij}$ is the symbol used in \cite{ng2008framed}, but not the $\infty \times \infty$ matrix we defined before. In the language of \cite{ng2008framed}, our $a_{ij}^0$ is the same as $a_{ij}'$ in that paper. \footnote{In the Ng's paper just mentioned, $a_{ij}'$ was defined differently. But we think that was an typo and our argument here is the right way to define it. Also note that the matrix $A$ in that paper has entries $a_{ij}'$, but not $a_{ij}$.} Moreover,
 
 $\Phim{\beta}(a_{i, n+1}^0)*a_{n+1,j}^{0} =\Phim{\beta}(a_{i, n+1})*a_{n+1,j}^{0}= \sum\limits ((\Phi_{\beta}^{L})_{ik}a_{k,n+1})*a_{n+1,j}^{0} = \sum\limits ((\Phi_{\beta}^{L})_{ik}a_{k,n+1}^0)*a_{n+1,j}^{0}  = \sum\limits  (\Phi_{\beta}^{L})_{ik}a_{k,j}^0$.
 
 Then we have 
 
 $A_{ij}^{00} - (\Lambda\Phiml{\beta}A)_{ij}^{00} = a_{ij}^0 - \lambda^{\delta_{i,1}}\Phim{\beta}(a_{i, n+1}^0)*a_{n+1,j}^{0} = a_{ij}^0 - \lambda^{\delta_{i,1}} \sum \limits (\Phi_{\beta}^{L})_{ik}a_{k,j}^0$,
 
 which is exactly the $(i,j)$-entry of $A - \Lambda \Phi_{\beta}^{L}A$ defined in \cite{ng2008framed}. Similarly, $A_{ij}^{00} - (A\Phimr{\beta}\Lambda^{-1})_{ij}^{00} = a_{ij}^0 - \lambda^{-\delta_{j,1}}a_{i, n+1}^0* \Phim{\beta}(a_{n+1,j}^{0})$ is the $(i,j)$-entry of $A - A\Phi_{\beta}^{R}\Lambda^{-1}$. Therefore, there is a well-defined isomorphism $H_0 \longrightarrow hc_0(\beta)$ sending $a_{ij}^0$ to $\mu a_{ij}$ if $i>j$ and $a_{ij}$ otherwise.

\end{proof}
\end{prop}

\begin{cor}
If $K$ is a local knot in $\Lens$ with framing $l$, then $HC_0(K;l)$ is infinitely generated as an $R$-algebra.
\begin{proof}
Clear from Proposition \ref{HC0 decomposition}.
\end{proof}
\end{cor}

We just showed that $HC_0$ is infinitely generated for local knots.  On the other hand, Theorem \ref{thm torus knot} shows $HC_0$ is always finitely generated for torus knots. Through some amount of computer calculations, we find that $HC_0$ is always finitely generated for non-local knots. This motivates us to come up with following conjecture.

\begin{conjecture} \label{conj finite generate}
Let $K$ be a knot in $\Lens$ with framing $l$, then $HC_0(K;l)$ is finitely generated as an $R$-algebra if and only if $K$ is not local.
\end{conjecture}

\subsection{Augmentations} \label{subsec:augmentation}

The invariant, $HC_0$, could be very difficult to compute for general knots, especially when the number of crossings is large. Thus we will deduce a family of invariants from $HC_0$, which are called augmentation numbers and which are relatively easier to compute, at least by computers. The concept of augmentation numbers are introduced in \cite{ng2005knotI} \cite{epstein2001chekanov} for basically the same reason.

Let $d \geq 2$ be an integer and let $\Z_d = \Z/d\Z$. Pick three invertible numbers $\lambda_0,\mu_0, \Var_0 \in \Z_d$. Then $\Z_d$ can be treated as an $R$-module, with $\lambda, \mu, \Var$ acting by multiplication by $\lambda_0,\mu_0, \Var_0$, respectively. Then $H(\beta;f;d;\lambda_0,\mu_0, \Var_0):= HC_0(\beta;f) \otimes _{R} \Z_d$ is a $\Z_d$-algebra. Assume $HC_0(\beta;f)$ is finitely generated, then $H(\beta;f;d;\lambda_0,\mu_0, \Var_0)$ is a finitely generated $\Z_d$-algebra, and thus has finitely many algebra morphisms into $\Z_d$.

\begin{definition}
Let $\beta \in C_n, f \in \Z, 2 \leq d \in \Z$ such that $HC_0(\beta;f)$ is finitely generated as an $R$-algebra, and let $\lambda_0,\mu_0, \Var_0 \in \Z_d$ be invertible, then $Aug(\beta;f;d;\lambda_0,\mu_0, \Var_0)$ is defined to be the number of algebra morphisms from $H(\beta;f;d;\lambda_0,\mu_0, \Var_0)$ to $\Z_d$.
\end{definition}

For example, denote the braid $(\alpha_0 \cdots \alpha_{p-1})^q$ representing the $(p,q)$-torus knot by $T(p,q)$, then $Aug(T(1,4);0;3;1,1,2) = 4, \\ Aug(T(1,5);0;3;1,1,2) = 2, \, Aug(T(1,6);0;3;1,1,2) = 4, \, \\  Aug(T(1,4);0;5;1,1,3) = 6, \, Aug(T(1,5);0;5;1,1,3) = 3.$


\begin{prop} \label{prop const map}
Set $\lambda = \mu = 1$, then for any $\beta \in \C_n$, there is a $\Z[\Var^{\pm 1}]$-algebra morphism from $HC_0(\beta;f)$ to $\Z[\Var^{\pm 1}]$ sending each $a_{ij}^x$ to $2\Var$.
\begin{proof}
Let $t: \A_n \longrightarrow R, \, t(a_{ij}^x) = 2\Var.$ We first show for $\beta \in \C_n$, $t \Phi_{\beta} = t$. Clearly, it suffices to prove $t \Phi_{\alpha_k} = t, 0 \leq k \leq n-1.$ This can be checked directly from Equations \ref{equ:Phi_k}, \ref{equ:Phi_0}.

Similarly, one can prove $t \Phip{\beta} = t \Phim{\beta} = t.$

We need to show $t$ factors through $\mathcal{I}_{\beta;f;1,1}$. Note that now $\Lambda_{\beta;f;1,1}$ is the identity matrix.

$t((\Phiml{\beta}A)_{ij}^{xy}) = t(\Phim{\beta}(a_{i,n+1}^x)*a_{n+1,j}^y) = t \Phim{\beta}(a_{i,n+1}^x) = t(a_{i,n+1}^x) = 2\Var = t(A_{ij}^{xy})$.

The other three relations can be verified analogously.
\end{proof}
\end{prop}

\begin{cor}
Let $\beta \in \C_n, f \in \Z$ let $\Var_0 \in \Z_d$ be invertible, then $Aug(\beta;f;d;1,1,\Var_0) \geq 1.$
\begin{proof}
The map $t$ defined in Proposition \ref{prop const map} naturally induces a map from $H(\beta;f;d;1,1,\Var_0)$ to $\Z_d$.
\end{proof}
\end{cor}

\section{A topological interpretation of the knot invariant} \label{sec:topology}

In this section, we show that the framed knot invariant $HC_0$ actually has a rather simple interpretation as the framed cord algebra given in Definition $2.2$ in \cite{ng2008framed}. The framed cord algebra is defined for an oriented framed knot $K$ in an oriented $3$-manifold $M$, which is conjectured to be the zero-th relative contact homology of $\Lambda_K$ in $ST^{*}M$. In the same paper, the author also gave a cord interpretation of the framed cord algebra for knots in $S^3$ with $0$ framing. In the following, we modify the cord interpretation so that it adapts to the knots with any framing, and prove that the modified version is equivalent to the framed cord algebra. Then we show that the knot invariant $HC_0$ coincides with the framed cord algebra.

\begin{definition} \label{cordint}
 Suppose $M$ is an oriented $3$-manifold, and $K$ an oriented framed knot in $M$ with $l,m$ the homotopy classes of the longitude and the meridian of $K$ in $\pi_1(M \setminus K)$. Fix a point $*$ on $l$.

1). A cord in $M$ relative to $(K,l)$ is a continuous map $\gamma: [0,1] \longrightarrow M \setminus K$, such that $\gamma(0), \gamma(1) \in l$ and $\gamma^{-1}(*) = \emptyset$. Two cords $\gamma_1, \gamma_2$ are said to be equivalent if they are homotopic relative to $l \setminus \{*\}$. Informally speaking, one can slide a cord $\gamma$ along $l$, so long as not to pass through the point $*$.

2). Let $R$ be the ring $\mathbb{Z}[\lambda^{\pm}, \mu^{\pm}, \Gamma^{\pm}]$.  The framed cord algebra, $A(K,l;M)$, is defined as the algebra over $R$ freely generated by the equivalent classes of cords, modulo the ideal generated by the relations given in Figure \ref{Skein}.

\begin{figure}[h!]
\centering
\def\svgwidth{11cm}
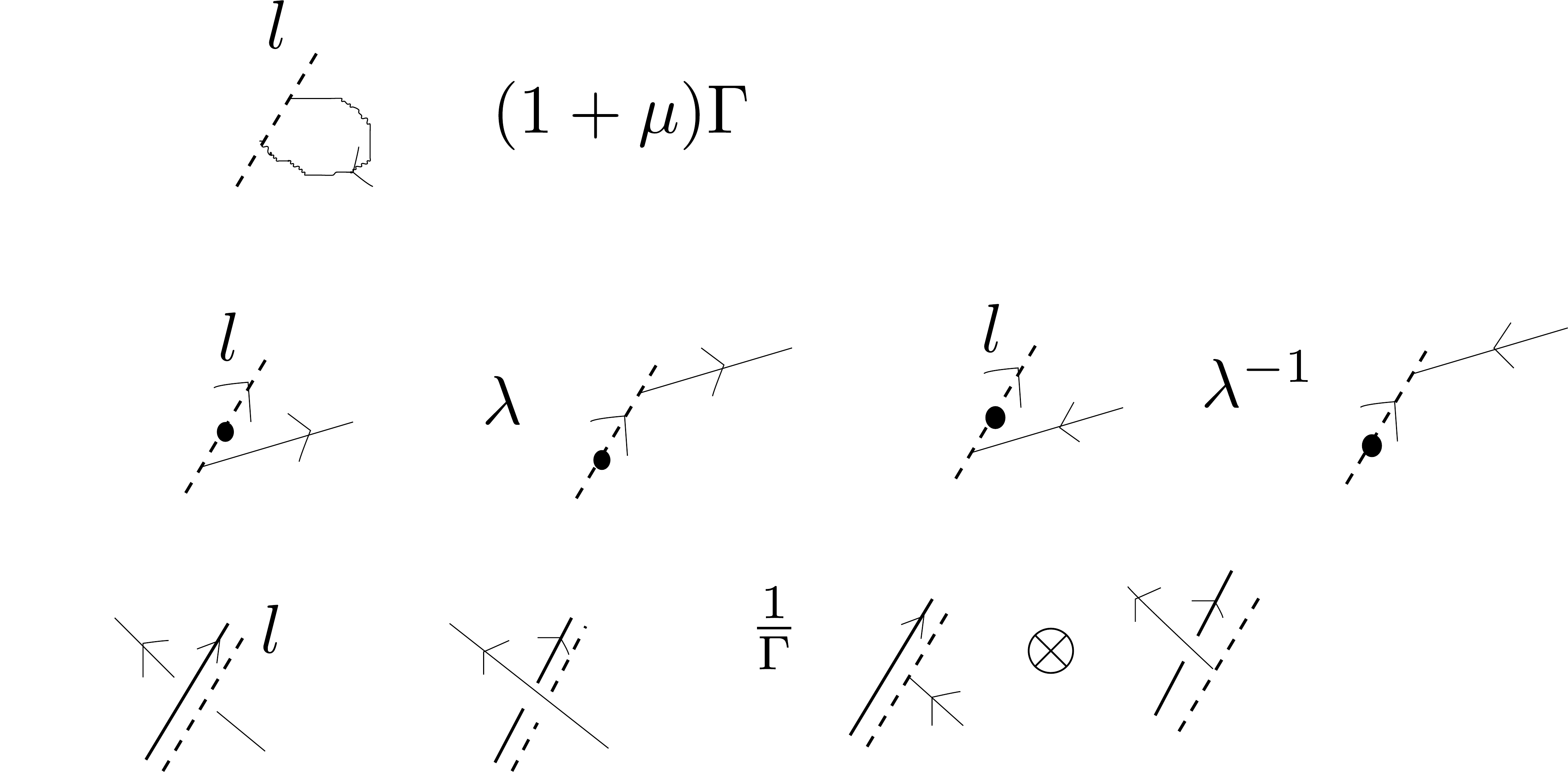
\caption{Skein relation}
 \label{Skein}
\end{figure}

\end{definition}

In Figure \ref{Skein}, the dashed line stands for the curve representing $l$, and the cord is represented by the solid line transversal to $l$ while the knot is drawn as the solid line parallel to $l$. In the third relation, the diagrams are understood to depict some local neighborhood outside of which the diagrams agree, and the meridian $m$ is assumed to rotate around $K$ counter clock-wise under the current projection.

Now we prove that the framed cod algebra is isomorphic to the one defined in \cite{ng2008framed}. For the readers convenience, we first recall the definition of framed cord algebra there.

\begin{definition}\cite{ng2008framed} \label{homopyint}
Let $K \subset M$ be an oriented framed knot in an oriented $3$-manifold $M$, and let $l,m$ denote the homotopy classes of the longitude and meridian of $K$ in $\pi_1(M \setminus K)$. The framed cord algebra, $\tilde{A}(K,l;M)$, of $K$ is the algebra over $R$ freely generated by the elements of $\pi_1(M \setminus K)$, modulo the ideal generated by the relations

$1). [e] = (1+\mu)\Var;$

$2). [\gamma l] = [l \gamma] = \lambda[\gamma]$ for $\gamma \in \pi_1(M \setminus K)$;

$3). [\gamma_1 \gamma_2 ] + [\gamma_1 m \gamma_2] = \frac{1}{\Var}[\gamma_1][\gamma_2],$ for $\gamma_1,\gamma_2 \in \pi_1(M \setminus K)$.
\end{definition}

 \begin{rem} \label{hotopycordrem}
 1). If we set $\gamma_1 = \gamma, \gamma_2 = e$, then from the first and the third relation, we can derive the relation $[\gamma m]  = \mu[\gamma].$ Similarly, we have $[m \gamma] = \mu[\gamma]$.

 2). If $[l'] = [l][m]^f$, then $\tilde{A}(K,l';M)$ can be obtained from $\tilde{A}(K,l;M)$ by replacing $\lambda$ by $\lambda\mu^{-f}$.
 \end{rem}

Clearly, the framed cord algebra does not depend on the choice of the base point in defining $\pi_1(M \setminus K)$.

\begin{prop}
The framed cord algebras defined in Definition \ref{cordint} and Definition \ref{homopyint} coincide, namely, $A(K,l;M) \simeq \tilde{A}(K,l;M)$ for an oriented knot $K$ with framing (longitude) given by $l$ in the manifold $M$.
\begin{proof}
Assume the base point $p$ is on the curve $l$, different from the point $*$. And for a point $z \in l$, let $\tau_z$ be the sub arc of $l$ connecting $p$ to $z$ not passing the point $*$.  Clearly, an element of $\pi_1(M \setminus K)$ is automatically an equivalence class of cords. Also it's easy to see that the three relations in defining $\tilde{A}(K,l;M)$ turn into the three relations defining $A(K,l;M)$, respectively. Conversely, for a cord $\gamma$, let $\tilde{\gamma} = \tau_{\gamma(0)} * \gamma * \bar{\tau}_{\gamma(1)}$. Then $\tilde{\gamma}$ is an element of $\pi_1(M \setminus K)$, and this map also preserves the defining relations.
\end{proof}
\end{prop}

\begin{thm}
Let $\beta \in \C_n$ be a braid whose closure is a knot in $\Lens$, and let $l,m$ be the homotopy classes of the longitude and the meridian of $\hat{\beta}$ in $\pi_1(\Lens \setminus \hat{\beta})$, such that $[l] = [\hat{\beta}'][m]^f$, where $\hat{\beta}'$ is a parallel copy diagram of $\hat{\beta},$ and $f \in \Z$ is an integer. Then we have $HC_0(\beta;f) \simeq A(\hat{\beta},l;\Lens)$.
\begin{proof}

By the second part of Remark \ref{hotopycordrem} and the properties of $HC_0(\beta;f)$, it suffices to prove the theorem for $f = 0$, namely $[l] = [\hat{\beta}']$. Set $\Lambda = \Lambda_{\beta;0;1,1}$.

Let $X = D_n \times [0,1] /\{ (x,0) \sim (x,1), x \in D_n\}$.  Present $\beta$ as a braid diagram inside $X$. See Figure \ref{ModelS1S2}. Assume $\beta$ intersect $D_n$ in $p_1, \cdots, p_n$. Take a parallel copy diagram $\beta'$ of $\beta$, such that $\beta'$ intersects with $D_n$ in the points $q_1, \cdots, q_n$. Also choose some point on $\hat{\beta'}$ right above $q_1$ as the point $*$.

It's clear that any cord in $\Lens$ relative to $(\hat{\beta},\hat{\beta'})$ can be homotoped to inside $X$. Then we slide the cord $\gamma$ along $\hat{\beta'}$ and whenever the cord passes the point $*$, we will multiply $\lambda(\lambda^{-1})$ to it according to the second relation in Figure \ref{Skein}. Finally the cord is slided into $D_n \times \{0\}$. we denote the resulting curve by $\tilde{\gamma}$, which is an element in $\Q_n$.

We define the map $\varphi: A(\hat{\beta},l;\Lens) \longrightarrow HC_0(\beta;0)$ by sending any cord $\gamma$ to $\lambda^{s}\tilde{\gamma}$, where $\lambda^{s}$ is the scalar gathered on the way to transit $\gamma$ into $\tilde{\gamma}$, as stated in the above paragraph. There are several points where we need to check the map is well-defined.

{\it{Step 1}}: The projection of $\gamma$ to $D_n \times \{0\}$ is not unique, and different projections differ by actions of $\Phi_{\beta}$. So we need to show for any $\gamma \in \Q_n$ from $q_i$ to $q_j$, we have $\gamma  = \lambda^{\delta_{i,1}}\Phi_{\beta}(\gamma)\lambda^{-\delta_{j,1}}$ in $HC_0(\beta;0)$. Since $\gamma$ can be written as a sum of monomials of the form $a_{i,i_1}^{x_1}a_{i_1,i_2}^{x_2} \cdots a_{i_{k-1},j}^{x_k}$,  by Corollary \ref{A-PhiA}, $\gamma  - \lambda^{\delta_{i,1}}\Phi_{\beta}(\gamma)\lambda^{-\delta_{j,1}}$ is contained in $\mathcal{I}_{\beta;0;1,1}$ and thus $0$ in $HC_0(\beta;0)$.

{\it{Step 2}}: In $\Lens$, the cords have more flexibilities to be homotoped than in $X$. Precisely, there are two more type of flexibilities. Let $\gamma_1, \gamma_2$ be two curves in $D_n$ such that $\gamma_1(1) = \gamma_2(0) = z_1, \gamma_1(0) = q_i, \gamma_2(1) = q_j, $ and let $\delta$ be the loop $\{z_1\} \times S^1$, then it's clear that $\gamma_1 * \gamma_2, \gamma_1 * \delta * \gamma_2$ are equivalent cords in $\Lens$ but not in $X$. If we project $\gamma_1 * \delta * \gamma_2$ to $D_n \times \{0\}$, then we get $\lambda^{\delta_{i,1}}\Phim{\beta}(\gamma_1)*\gamma_2 $ or $\gamma_1 * \Phim{\beta}(\gamma_2)\lambda^{-\delta_{j,1}}$. These are guaranteed by the relations $A - \Lambda\Phiml{\beta} A, A - A\Phimr{\beta}\Lambda^{-1}$. See Part $(2)$ of Remark \ref{rem relation}.

Similarly, in the above argument, if we replace \lq\lq $z_1$" by \lq\lq $z_0$", then we get the relations $A - \Lambda\Phipl{\beta} A, A - A\Phipr{\beta}\Lambda^{-1}$.

{\it{Step 3}}: The first and the third relation in Figure \ref{Skein} that defines $A(\hat{\beta},l;\Lens)$ are apparently mapped to the two \lq\lq skein" relations that define $\tilde{A}_n$. And the second relation in the same figure is also preserved by the map.

The above three steps showed that $\varphi$ is well-defined. It's also easy to prove it's a bijection.

\end{proof}
\end{thm}

\bibliographystyle{plain}
\bibliography{KCHbib}

\begin{thebibliography}{10}

\bibitem{Ding}
Feifei Chen, Fan Ding, and Youlin Li.
\newblock Legendrian torus knots in ${S}^1\times {S}^2$.
\newblock 2013.
\newblock arxiv.org/abs/1310.1535.

\bibitem{crisp1999injective}
John Crisp.
\newblock Injective maps between {Artin} groups.
\newblock In {\em Geometric Group Theory Down Under: Proceedings of a Special
  Year in Geometric Group Theory, Canberra, Australia, 1996}, page 119. Walter
  de Gruyter, 1999.

\bibitem{epstein2001chekanov}
Judith Epstein, Dmitry Fuchs, and Maike Meyer.
\newblock Chekanov--{E}liashberg invariants and transverse approximations of
  {L}egendrian knots.
\newblock {\em Pacific J. Math}, 201(1):89--106, 2001.

\bibitem{HeltonPackage}
J~Helton, R~Miller, and M~Stankus.
\newblock {NCAlgebra}: a {Mathematica} package for doing non-commuting algebra.
\newblock \url{http://www.math.ucsd.edu/~ncalg}.

\bibitem{Lin}
Xiaosong Lin.
\newblock Markov theorems for links in 3-manifolds.
\newblock {\em Physics and Topology, Nankai Tracts in Mathematics}, 12:360,
  2007.

\bibitem{ng2005knotI}
Lenhard Ng.
\newblock Knot and braid invariants from contact homology {I}.
\newblock {\em Geom. Topol}, 9:247--297, 2005.

\bibitem{ng2005knot}
Lenhard Ng.
\newblock Knot and braid invariants from contact homology {II}.
\newblock {\em Geom. Topol}, 9:1603--1637, 2005.

\bibitem{ng2008framed}
Lenhard Ng.
\newblock Framed knot contact homology.
\newblock {\em Duke Mathematical Journal}, 141(2):365--406, 2008.

\bibitem{ng2011combinatorial}
Lenhard Ng.
\newblock Combinatorial knot contact homology and transverse knots.
\newblock {\em Advances in Mathematics}, 227(6):2189--2219, 2011.

\bibitem{CuiPackage}
Shawn X.~Cui.
\newblock A computer package computing the ${HC}_0$ invariant and related
  invariants for knots in ${S}^1 \times {S}^2$.
\newblock 2014.
\newblock \url{http://math.ucsb.edu/~xingshan/publication.html}.

\end{thebibliography}

\end{document}